\documentclass[12pt]{article}
\usepackage{amsthm,amsmath, amssymb, amsfonts, amscd, xspace, pifont}
\usepackage{epsfig, psfrag, psfig, pstricks}
\usepackage{float, fancybox, fullpage}
\usepackage{delarray}
\usepackage{hyperref}
\usepackage{url}
\usepackage{anysize}
\usepackage{multirow}
\usepackage[top=1.5in, bottom=2in, left=0.95in, right=0.95in]{geometry}
\usepackage{rotating}
\usepackage{latexsym,ifthen}
\usepackage{subfigure}
\usepackage{anysize}
\usepackage{natbib}
\usepackage{times, url}
\usepackage{bm}
\usepackage{natbib,amsmath,epsfig,epsf,psfrag,  amssymb, amsfonts}

\begin{document}

\newtheorem{remark}{\indent \sc Remark}
\newtheorem{theorem}{\indent \sc Theorem}
\newtheorem{lemma}{\indent \sc Lemma}

\newcommand {\apgt} {\ {\raise-.5ex\hbox{$\buildrel>\over\sim$}}\ }
\newcommand {\aplt} {\ {\raise-.5ex\hbox{$\buildrel<\over\sim$}}\ }
\def\bv{\mathbf b}
\newcommand{\phif}{\textsc{igf}}
\newcommand{\sign}{\mbox{sign}}
\newcommand{\phgf}{\textsc{gf}}
\newcommand{\fix}{$\textsc{gf}_0$}
\newcommand{\mb}[1]{\mbox{\bf #1}}
\newcommand{\Exp}[1]{\mbox{E}\left[#1\right]}
\newcommand{\pr}[1]{\mbox{P}\left[#1\right]}
\newcommand{\pp}[0]{\mathbb{P}}
\newcommand{\ee}[0]{\mbox{E}}
\newcommand{\re}[0]{\mathbb{R}}
\newcommand{\argmax}[0]{\mbox{argmax}}
\newcommand{\argmin}[0]{\mbox{argmin}}
\newcommand{\ind}[0]{\mbox{\Large\bf 1}}
\newcommand{\narrow}{\stackrel{n\rightarrow\infty}{\longrightarrow}
}
\newcommand{\weakpn}{\stackrel{P_n}{\leadsto}}
\newcommand{\weakpnboot}{\mbox{\raisebox{-1.5ex}{$\stackrel
{\mbox{\scriptsize $P_n$}}{\stackrel{\mbox{\normalsize
$\leadsto$}} {\mbox{\normalsize $\circ$}}}$}}\,}
\newcommand{\ol}[1]{\overline{#1}}
\newcommand{\avgse}[1] { \bar{\widehat{\sigma}}_{#1} }
\newcommand{\mcse}[1]  { \sigma^{*}_{#1} }
\newcommand{\po}{\textsc{po}}
\newcommand{\ph}{\textsc{ph}}



\title{Semiparametric Additive
Transformation Model under Current Status Data}

\author{Guang Cheng\thanks{Corresponding Author, Purdue University, West Lafayette, IN 47907,
Email: chengg@purdue.edu.}\; and Xiao Wang\thanks{Purdue University, West Lafayette, IN 47907,
Email: wangxiao@purdue.edu.}}
\date{Purdue University}

\maketitle

\begin{abstract}
We consider the efficient estimation of the semiparametric additive transformation model with current status data. A wide range of survival models and econometric models can be incorporated into this general transformation framework. We apply the B-spline approach to simultaneously estimate the linear regression vector, the nondecreasing transformation function, and a set of nonparametric regression functions. We show that the parametric estimate is semiparametric efficient in the presence of multiple nonparametric nuisance functions. An explicit consistent B-spline estimate of the asymptotic variance is also provided. All nonparametric estimates are smooth, and shown to be uniformly consistent and have faster than cubic rate of convergence. Interestingly, we observe the convergence rate interfere phenomenon, i.e., the convergence rates of B-spline estimators are all slowed down to equal the slowest one. The constrained optimization is not required in our implementation. Numerical results are used to illustrate the
finite sample performance of the proposed estimators.
\end{abstract}

\vspace{0.01in} \noindent {\em Key Words:} B-spline; Consistent variance estimation; Current status data; Efficient estimation;
Semiparametric transformation models

\section{Introduction}

We consider the efficient estimation of the following semiparametric additive transformation model:
\begin{eqnarray}
H(U)=Z'\beta+\sum_{j=1}^{d}h_{j}(W_{j})+\epsilon,\label{trans}
\end{eqnarray}
where $H(\cdot)$ is a monotone transformation function, $h_{j}(\cdot)$'s are smooth regression functions (with possibly different degrees of smoothness), and
$\epsilon$ has a known distribution $F(\cdot)$ with support $\mathbb{R}$. A wide range of survival models and econometric models can be incorporated into the above general transformation framework, e.g., \citep{hr97, s98, h99, bbg06, bmm08}. In particular, the model (\ref{trans}) can be readily applied to a failure time $T$ by letting $U=\log T$. We can obtain the partly linear additive Cox model, i.e., \cite{h99}, by assuming $F(s)=1-\exp(-e^{s})$ and $H(u)=\log A(e^{u})$,
where $A$ is an unspecified cumulative hazard function. Specifically, the hazard function of $T$, given the covariates $(z,w)$, has the form
\begin{eqnarray}
\lambda(t|z,w)=a(t)\exp(\tilde{\beta}'z+\sum_{j=1}^{d}
\tilde{h}_{j}(w_{j})),\label{coxhaz}
\end{eqnarray}
where $a(t)$ is the baseline hazard function,
$\tilde{\beta}=-\beta$ and $\tilde{h}_{j}=-h_{j}$. However, if we change the form of $F(s)$ to $e^{s}/(1+e^{s})$, the model (\ref{trans}) just becomes the partly linear additive proportional odds model.

Motivated by the close connection with survival models, we focus on the current status data in this paper which arises not only in survival analysis but also in demography, epidemiology, econometrics and bioassay. More specifically, we observe $X=(V,\Delta,Z,W)$, where $V\in\mathbb{R}$ is a random examination time and $\Delta=1\{U\leq V\}$. We assume that $U$ and $V$ are independent given $(Z,W)$. Under current status data, the model (\ref{trans}) is also related to the semiparametric binary model studied in econometrics. Using the link function $F(\cdot)$, we assume that the probability of $\Delta=1$, given the covariates $(Z,W,V)$, is of
the expression:
\begin{eqnarray}
P(\Delta=1|Z,W,V)=F\left(\tilde{\beta}'Z+\sum_{j=1}^{d}\tilde{h}_{j}(W_{j})+H(V)\right).\label{link}
\end{eqnarray}
Note that \cite{bbg06} and \cite{bmm08} have done a great
deal of statistical estimation and hypothesis testing on the model (\ref{link}) (without $\widetilde h_{j}$ terms) by assuming $F(\cdot)$ to be
log-log function and logistic function, respectively. An extensive discussions on the relation between (\ref{link}) and survival models
can be found in \cite{dg90}. Recently a similar transformation model has been considered by \cite{ct10} but for the {\it right censored data}. They showed that the monotone transformation function is root-n estimable which will never be achieved in the case of current status data. This is the key theoretical difference between the two types of survival data.

In this paper, we employ the B-spline approach to simultaneously estimate the vector $\beta$, monotone $H$ and smooth $h_j$'s. The corresponding estimates are denoted as $\widehat\beta$, $\widehat H$ and $\widehat h_j$. In contrast, \cite{mk05b} apply the penalized NPMLE approach to (\ref{trans}) (with $d=1$) which yields a non-smooth step function $\check H$ and the penalized estimate $\check h$. Our B-spline framework has the following theoretical and computational advantages over the existing penalized NPMLE approach:
\begin{enumerate}
\item Our B-spline estimate $\widehat H$ is smooth and uniformly consistent. However, $\check H$ is always discontinues (regardless of the smoothness of its true function $H_0$) and has a bias which does not vanish asymptotically. More importantly, the convergence rate of our $\widehat H$ $(\widehat h)$ is shown to be faster than that of $\check H$ $(\check h)$, i.e., $O_P(n^{-1/3})$. Therefore, we expect more accurate inferences drawn from $\widehat H$ $(\widehat h)$.

\item We are able to give an explicit B-spline estimate for the asymptotic covariance of $\widehat\beta$ based on which the asymptotic confidence interval of $\beta$ can be easily constructed. Under very weak conditions, its consistency is proven. However, the block jackknife approach in Ma \& Kosorok (2005) requires more computation, and is even not theoretically justified.

\item Our spline estimation algorithm requires much less computation than the isotonic type algorithm used in Ma \& Kosorok (2005) since the order of jumps in the step function is supposed to be much larger than the order of knots we choose for estimating $H$ and $h_j$'s.

\end{enumerate}
Despite the non-root-n convergence rates of $\widehat H$ and $\widehat h_j$'s, we are able to show that $\widehat\beta$ is root-n consistent, asymptotically normal and semiparametric efficient. We derive the efficient information bound  by taking the general two-stage projection approach from Sasieni (1992) which is needed due to the involvement of multiple nonparametric functions in semiparametric models.
Interestingly, we observe the convergence rate interfere phenomenon for the B-spline estimators, i.e., the convergence rates of nonparametric estimators are all slowed down to equal the slowest one. Moreover, by approximating $\log\dot H$ with the B-spline, we can avoid the monotonicity constraint in the implementation, which is usually required in the literature, e.g., \cite{zhh10}.

The remainder of the paper is organized as follows. Section~\ref{estpro} describes the B-spline estimation procedure. The asymptotic properties such as
consistency and convergence rates of the estimates are obtained in Section~\ref{asyrate}. The asymptotic distribution of the parametric component is studied in Section~\ref{weakcon},
and its efficient information and the corresponding explicit B-spline estimate are given in Section~\ref{varest}. Simulation studies are presented in Section~\ref{simulat}.
We close with an appendix containing technical details.

\section{Semiparametric B-spline Estimation}\label{estpro}

\subsection{Assumptions}
We first define some notations. For any vector $v$, $v^{\otimes 2}=vv'$. The notations $\apgt$ and $\aplt$ mean
greater than, or smaller than, up to a universal constant. We denote $A_n\asymp B_n$ if $A_n\aplt B_n$ and $A_n\apgt B_n$. The
notations $\pp_n$ and $\mathbb{G}_n$ are used for the empirical distribution and the empirical process of the observations,
respectively. Furthermore, we use the operator notation for evaluating expectation. Thus, for every measurable function $f$
and true probability $P$,
\begin{displaymath}
\pp_{n}f=\frac{1}{n}\sum_{i=1}^{n}f(X_i), \;\;\; Pf=\int{f}dP
\;\;\; \mbox{and} \;\;\;
\mathbb{G}_{n}f=\frac{1}{\sqrt{n}}\sum_{i=1}^{n}{(f(X_i)-Pf)}.
\end{displaymath}

We next present some model assumptions.
\begin{enumerate}
\item[M1.] $U$ and $V$ are independent given $(Z,W)$.

\item[M2.] (a) The covariates $(Z, W)$ are assumed to belong to a
bounded subset in $\mathbb{R}^{l+d}$, say
$[0,1]^{l}\times[0,1]^{d}$. The support for $V$ is
$[l_{v},u_{v}]$, where $-\infty<l_{v}<u_{v}<+\infty$; (b) The joint density for $(Z,V,W)$  w.r.t. Lebesgue measure stays away from zero, and the joint density for $(V,W)$ stays away from infinity.

\item[M3.] $E(Z-E(Z|V,W))^{\otimes 2}$ is strictly positive definite.

\item[M4.] The residual error distribution $F(\cdot)$ is assumed
to be known and has support $\mathbb{R}$. Denote the first,
second and third derivative of $F$ as $f$, $\dot{f}$ and $\ddot f$, respectively. We
assume that (a) $(f(u)\vee|\dot f(u)|\vee|\ddot f(u)|)\leq M<\infty$ over the whole $\mathbb{R}$ and $f(u)$ stays away from zero in any
compact set of $\mathbb{R}$; (b) $[f^{2}(v)-\dot{f}(v)F(v)]\wedge[f^{2}(v)+\dot{f}(v)(1-F(v))]>0,$
for all $v\in\mathbb{R}$.
\end{enumerate}
Since we employ the smooth B-spline estimation rather than the penalized NPML estimation, our residue error Condition M4 is much less restrictive than that in \cite{mk05b}, and may apply to more general class of semiparametric transformation models. Note that Condition M4(b) ensures the concavity of the function $s\mapsto \delta\log F(s)+(1-\delta)\log (1-F(s))$ for $\delta=0,1$.

It is easy to verify that the above Condition M4 is satisfied in the following two general classes of residue error distribution
functions after some algebra.
\begin{enumerate}
\item[F1.]
$F(s)=\gamma[2\Gamma(\gamma^{-1})]^{-1}\int^{s}_{-\infty}\exp(-|t|^
{\gamma})dt$ for $\gamma>1$ is a family of distributions, which includes the standard normal distribution after appropriate rescaling ($\gamma=2$). This corresponds to the probit model \cite{kp80}.

\item[F2.]  $F(s)=1-[1+\gamma e^{s}]^{-1/\gamma}$ is a Pareto
distribution with parameter $\gamma\in(0,\infty)$ and corresponds
to the odds-rate transformation family, see \cite{dd88a, dd88b}.
It includes the following two well-known special cases:
\begin{enumerate}
\item[(a).] Given $\gamma\rightarrow 0$, it yields the extreme
value distribution, i.e. $F(s)=1-\exp(-e^{s})$, which corresponds
to the complementary log-log transformation, see \cite{bbg06};

\item[(b).] Given $\gamma=1$, it gives the logistic distribution,
i.e. $F(s)=e^{s}/(1+e^{s})$, which corresponds to the logit
transformation, see \cite{bmm08}.
\end{enumerate}
\end{enumerate}

\subsection{B-spline Estimation Framework}
From now on, we change the signs of $\beta$ and $h_{j}$ for
simplicity of exposition. In addition, we re-center $H(v)$ to $H(v)-H(l_v)$ so that $H(l_v)=0$ for the purpose of identifiability. The additional parameter $H(l_v)$ will be absorbed into the vector $\beta$, i.e., the first coordinate of $z$ is set as one. Given a single observation at $x=(v,\delta,z,w)$, the log-likelihood of model (\ref{trans}) is written as
\begin{eqnarray}
\ell(\beta,h_1,\ldots,h_d,H)&=&\delta\log\left\{F\left[H(v)+\beta'z+\sum_{j=1}^{d}h_{j}(w_{j})\right]\right\}\nonumber
\\&&+(1-\delta)\log\left\{1-
F\left[H(v)+\beta'z+\sum_{j=1}^{d}h_{j}(w_{j})\right]\right\}.\label{log-lik}
\end{eqnarray}

We assume that $\beta\in\mathcal{B}$, which is a bounded
open subset in $\mathbb{R}^{l}$, and that its true value $\beta_0$ is an interior point of $\mathcal{B}$. Before specifying the parameter spaces for $H$ and $h_j$'s, we first introduce the H\"{o}lder ball $\mathbf{H}^r_c(\mathcal{Y})$, which is a class of smooth functions
widely used in the nonparametric
estimation, e.g., \cite{s82, s85}. For any $f\in\mathbf{H}^r_c(\mathcal{Y})$, it is $J<r$ times
 continuously differentiable on $%
\mathcal{Y}$ and its $J$-th derivative is uniformly H\"{o}lder continuous with exponent $\kappa \equiv r-J\in (0,1]$, i.e., $$\sup_{y_1,y_2\in\mathcal{Y}, y_1\neq y_2}
\frac{|f^{(J)}(y_1)-f^{(J)}(y_2)|}{|y_1-y_2|^\kappa}\leq c.$$
The functions in the H\"{o}lder ball can always be approximated by a basis expansion, i.e.,
\begin{eqnarray}
f(t)\approx\sum_{k=1}^{K}\gamma_k B_k(t)=\gamma'\mathbf{B}(t),\label{bsplapp}
\end{eqnarray}
where $\gamma=(\gamma_{1},\ldots,\gamma_{K})'$ and $\mathbf{B}(t)=(B_1(t),\ldots,B_{K}(t))'$. Actually, if the degree $d$ of the B-spline satisfies $d\geq(r-1)$, we have
\begin{eqnarray}
\|f-\gamma'\mathbf{B}\|_\infty\asymp K^{-r}\;\;\;\;\mbox{as}\;K\rightarrow\infty,\label{sieapp}
\end{eqnarray}
where $\|\cdot\|_{\infty}$ denotes the supremum norm..

Assume the following parameter space Condition P1 for the smooth $h_j$.
\begin{enumerate}
\item[P1.] For $j=1,\ldots,d$ and some known $c_j$, we assume that the parameter space for $h_j$ is $\mathcal{H}_j$, where
$$\mathcal{H}_j=\left\{h_j: h_j\in\mathbf{H}^{r_j}_{c_j}[0,1]\;\mbox{with}\;r_j>1/2\;\mbox{and}\;\int_0^1 h_j(w_j)dw_j=0\right\},$$ and that the corresponding spline space is
$$\mathcal{H}_{jn}=\left\{h_j: h_j(w)=\gamma_j'\mathbf{B}_j(w)\;\mbox{with}\;\|h_j\|_\infty\leq c_j\;\mbox{and}\;\int_0^1 h_j(w_j)dw_j=0\right\},$$ based on a system of basis functions $\mathbf{B}_j=(B_{j1},\ldots,B_{jK_j})'$ of degree $d_j\geq (r_j-1)$.
\end{enumerate}

As seen from the previous examples, it is reasonable to assume that $H(\cdot)$ is differentiable and strictly increasing over $[l_v,u_v]$, i.e., $\dot H(v)\geq C_0>0$.
Considering that $H(l_v)=0$,  we can thus write $H(v)=\int_{l_v}^v\exp(g(s))ds$, where $g(v)\equiv\log \dot H(v)$ is well defined.
Such reparametrization can get around the strict monotonicity and positivity constraints of $H$, and thus avoids the constrained optimization
in the computation. The parameter space Condition P2 for $g$ is specified below.
\begin{enumerate}
\item[P2.] For some known $c_0$, we assume that the parameter space for $g$ is $\mathcal{G}$,
where $$\mathcal{G}=\left\{g: g\in
\mathbf{H}^{r_0}_{c_0}[l_v,u_v]\;\mbox{with}\;r_0>1/2\right\},$$ and that the corresponding spline space is $$\mathcal{G}_n=\left\{g: g(v)=\gamma_0'\mathbf{B}_0(v)\;\mbox{and}\;\|g\|_\infty\leq c_0\right\}$$ based on a system of basis functions $\mathbf{B}_0=(B_{01},\ldots,B_{0K_0})$ of degree $d_0\geq (r_0-1)$.
\end{enumerate}
Similarly, we define $\mathcal{G}_n'=\{H(v)=\int_{l_v}^v\exp(g(s))ds: g\in\mathcal{G}_n\}$. By some algebra, we can show that $H\in\mathbf{H}_{c_0'}^{r_0+1}[l_v,u_v]$ for some $c_0'<\infty$.

\begin{remark}
Note that in the theoretical proofs and numerical calculations the exact values of $c_j$ are not necessary. Instead, only the boundedness condition,
equivalently the compactness of parameter spaces and spline spaces, is needed. Here we assume this boundedness condition, which can be relaxed by invoking
the chaining arguments,
only for simplifying our theoretical derivations.
\end{remark}

In this paper, we propose the B-spline approach to estimate $H$ and $h_j$'s as follows. Let $\mathcal{A}=\mathcal{B}\times\mathcal{G}\times\Pi_{j=1}^{d}\mathcal{H}_j$ and $\mathcal{A}_n=\mathcal{B}\times\mathcal{G}_n\times\Pi_{j=1}^{d}\mathcal{H}_{jn}$. Denote $\alpha$ as $(\beta',g,h_1,\ldots,h_d)'$ and its true value $\alpha_0$ as $(\beta_0', g_0, h_{10},\ldots,h_{d0})'$, where $g_0(\cdot)=\log \dot H_0(\cdot)$. The log-likelihood (\ref{log-lik}) for the observation $i$ can thus be reparametrized as
\begin{eqnarray}
\ell_i(\alpha)&=&\delta_i\log\left\{F\left[\beta'z_i+\int_{l_v}^{v_i}\exp(g(s))ds+\sum_{j=1}^{d}h_j(w_{ij})\right]\right\}\nonumber\\
&&+(1-\delta_i)\log\left\{1-F\left[\beta'z_i+\int_{l_v}^{v_i}\exp(g(s))ds+\sum_{j=1}^{d}h_j(w_{ij})\right]\right\}.\label{sie-lik}
\end{eqnarray}
The corresponding B-spline estimate $\widehat\alpha$ is defined as
\begin{eqnarray}
\widehat\alpha=\arg\max_{\alpha\in\mathcal{A}_n}\sum_{i=1}^n\ell_i(\alpha).\label{sieest}
\end{eqnarray}
We can also write $\widehat\alpha=(\widehat\beta',\widehat g,\widehat h_1,\ldots,\widehat h_d)'=(\widehat\beta',\widehat\gamma_0'\mathbf{B}_0,\widehat\gamma_1'\mathbf{B}_1,\ldots,\widehat\gamma_d'\mathbf{B}_d)'$. Then, the estimate  $\widehat H(v)=\int_{l_v}^v\exp(\widehat\gamma_0'\mathbf{B}_0(s))ds$. Some tedious algebra reveals that the Hessian matrix of $\ell_i(\alpha)$
w.r.t. $(\beta',\gamma_0',\gamma_1',\ldots,\gamma_d')'$ is indeed negative semidefinite under Condition M4(b) which guarantees the existence of $\widehat\alpha$. See more discussions on the computation feasibility in the simulation section. The above estimation procedure also applies to other linear sieves approximating the H\"{o}lder ball (or more generally H\"{o}lder space), e.g., wavelets.

\section{Consistency and Rates of Convergence}\label{asyrate}
In this section, we show that our B-spline estimate is consistent and the convergence rate of each nonparametric estimate appears to interfere with
 each other.
Define $$d(\alpha,\alpha_0)=\|\beta-\beta_0\|+\|H-H_0\|_2+\sum_{j=1}^{d}\|h_j-h_{j0}\|_2,$$ where $\|\cdot\|_2$ is the $L_2$ norm. Now we give the main Theorem of this section.
\begin{theorem}\label{consis}
Suppose that Conditions M1-M4 and P1-P2 hold. If $K_j/n\rightarrow 0$ for $j=0,1,\ldots,d$, then we have
\begin{eqnarray}
d(\widehat\alpha, \alpha_0)=o_P(1).\label{consfor}
\end{eqnarray}
More specifically, we further prove that
\begin{eqnarray}
d(\widehat\alpha,\alpha_0)=O_P\left(\max_{0\leq j\leq d}\left\{K_j^{-r_j}\vee\sqrt{K_j/n}\right\}\right).\label{convfor0}
\end{eqnarray}
If we further require that $K_j\asymp n^{1/(2r_j+1)}$ for $j=0,\ldots,d$, then we have
\begin{eqnarray}
d(\widehat\alpha,\alpha_0)=O_P(n^{-r/(2r+1)}),\label{convfor}
\end{eqnarray}
where $r=\min_{0\leq j\leq d}\{r_j\}$.
\end{theorem}
According to Theorem~\ref{consis}, the smooth $\widehat H$ can achieve the faster convergence rate, i.e., $O_P(n^{-r/(2r+1)})$, than $n^{1/3}$-rate
derived in the penalized estimation context, see \cite{mk05b}, when we assume that $g_0$ and $h_{j0}$'s are all at least continuously differentiable,
i.e., $r>1$. More importantly, we can further show that $\widehat H$ is uniformly consistent, i.e., $\|\widehat H-H_0\|_\infty=o_P(1)$, by applying Lemma 2 in \cite{cs98} that
 $\|f\|_\infty\aplt\|f\|_{L_2(Leb)}^{2r/(2r+d)}$ for any $f\in\mathbf{H}_c^r[a,b]^d$ and noting that  $\widehat H, H_0\in\mathbf{H}_{c_0'}^{r_0+1}[l_v,u_v]$ for some $c_0'>0$.

The above theorem also holds when we  employ the constrained monotone B-spline to approximate $H_0$, i.e., $\gamma_0'\mathbf{B}_0(v)\approx\log H(v)$
with $\gamma_{01}\leq\gamma_{02}\leq\ldots\leq\gamma_{0K_0}$. However, such constrained optimization usually requires additional computational effort,
see Zhang et al. (2010).

\begin{remark}
From the above Theorem~\ref{consis}, we observe the interesting convergence rate interfere phenomenon, i.e., the convergence rate for each
B-spline estimate is forced to equal the slowest one. In \cite{mk05b}, they also show that the convergence rate of
the penalized estimate $\widetilde h$ is unfortunately slowed down to $O_P(n^{-1/3})$ by the NPMLE $\widetilde H$ regardless of the smoothness degree of $h_0$.
One possible solution in achieving the optimal rate for
each nonparametric estimate is to extend the most recent mixed rate asymptotic results \cite{r08} to the semiparametric setup.
\end{remark}
Since we assume that $r>1/2$, the convergence rate given in (\ref{convfor}) is always $o_P(n^{-1/4})$.
Such a rate is usually fast enough to guarantee the regular asymptotic behavior of $\widehat\beta$, i.e., $\sqrt{n}$-consistency and asymptotic normality.
Indeed, we will improve the current suboptimal rate of $\widehat\beta$ in (\ref{convfor}) to the optimal $\sqrt{n}$ rate, and further show that $\widehat\beta$ is
semiparametric efficient in next section.

\section{Weak Convergence of the Parametric Estimate}\label{weakcon}
In this section, we study the weak convergence of the spline estimate $\widehat\beta$ in the presence of multiple nonparametric nuisance functions.
We first calculate the semiparametric efficient information based on the projection onto the nonorthogonal sumspace.

Let
\begin{eqnarray*}
Q_{\theta}(x)=f(\theta)\left(\frac{\delta}{F(\theta)}-\frac
{1-\delta}{1-F(\theta)}\right),
\end{eqnarray*}
where $\theta(z,v,w)=\beta'z+H(v)+\sum_{j=1}^d h_j(w_j)$. Denote $\theta_0$ as the true value of $\theta$.
The score functions (operators) for $\beta$, $g$ and $h_j$ are separately calculated as
\begin{eqnarray}
\dot{\ell}_{\beta}(X;\alpha)&=&ZQ_{\theta}(X),\label{scobeta}\\
\dot\ell_g[a](X;\alpha)&=&\left[\int_{l_v}^{V}\exp(g(s))a(s)ds\right]Q_{\theta}(X),\label{scoh}\\
\dot\ell_{h_j}[b_j](X;\alpha)&=&b_{j}(W_{j})Q_{\theta}(X).\label{scohj}
\end{eqnarray}
We assume that $a\in L_2(H)\equiv\{a: \int_{l_v}^{u_v}a^2(s)dH(s)<\infty\}$ and $b_j\in L_2^0(w_j)\equiv\{b_j:\int_{0}^{1}b_j(w_j)dw_j=0\;
\mbox{and}\;\int_{0}^{1}b_j^2(w_j)dw_j<\infty\}$ so that all the score functions defined above are square integrable.

To calculate the
efficient score function $\widetilde \ell_\beta$, we need to find the projection of $\dot\ell_\beta$ onto the sumspace $\mathbf{A}=A_g+A_{h_1}+
\cdots+A_{h_d}$, where
$A_{g}=\{\dot\ell_g[a]: a\in L_2(H)\}$ and $A_{h_j}=\{\dot\ell_{h_j}[b_j]: b_j\in L_2^0(w_j)\}$. For simplicity, we define $\dot\ell_{\beta}(X;\alpha_0)$ and $\dot\ell_\beta(X;\widehat\alpha)$ as $\dot\ell_{\beta_0}$ and
$\dot\ell_{\widehat\beta}$, respectively. The same notation rule applies to $\dot\ell_{g}[a](X;\alpha)$ and $\dot\ell_{h_j}[b_j](X;\alpha)$.
We define $$\widetilde\ell_{\beta}(X;\alpha)=
\dot\ell_{\beta}(X;\alpha)-\dot\ell_{g}[\bar a^\dag](X;\alpha)-\sum_{j=1}^d\dot\ell_{h_{j}}[\bar b_j^\dag](X;\alpha),$$ where $\bar a^\dag=
(a_1^\dag,\ldots,a_l^\dag)'$ and $\bar b_j^\dag=(b_{j1}^\dag,\ldots,b_{jl}^\dag)'$. And $(a_k^\dag,b_{1k}^\dag,\ldots,b^\dag_{dk})$ is the
minimizer of
$$(a_k,b_{1k},\ldots,b_{dk})\mapsto E\left\{[\dot\ell_{\beta_0}]_k-\dot\ell_{g_0}[a_k]-\sum_{j=1}^{d}\dot\ell_{h_{j0}}[b_{jk}]\right\}^2$$ for $k=1,\ldots,l$. Similarly, denote $\widetilde\ell_\beta(X;\alpha_0)$ and $\widetilde\ell_{\beta}(X;\widehat\alpha)$ as $\widetilde\ell_{\beta_0}$ and
$\widetilde\ell_{\widehat\beta}$, respectively. By taking the two-stage projection approach from Sasieni (1992), we have
\begin{eqnarray}
\tilde{\ell}_{\beta_{0}}(X)=
\left(Z-\bar b^{\dag}(W)-\frac{E((Z-\bar b^{\dag}(W))Q_{\theta_{0}}^{2}(X)
|V)}{E(Q_{\theta_{0}}^{2}(X)|V)}\right)Q_{\theta_{0}}(X)\label{effsco}
\end{eqnarray}
where $\bar b^{\dag}(W)=\sum_{j=1}^{d}\bar b_{j}^{\dag}(W_{j})$ satisfies
\begin{eqnarray}
E\left\{\left[Z-\bar b^{\dag}(W)-\frac{E((Z-\bar b^\dag(W))Q_{\theta_0}^2|V)}{E(Q_{\theta_0}^2|V)}\right]_kQ_{\theta_{0}}^{2}b_{jk}(W_j)\right\}
=0\label{inter0}
\end{eqnarray}
for every $b_{jk}\in L_2^0(w_j)$, $j=1,\ldots,d$ and $k=1,\ldots,l$. By slightly modifying the proof of Lemma 4 in \cite{mk05b}, we can show
that the above nonorthogonal projection is well defined and $\bar b^\dag(\cdot)$ exists by the alternating projection Theorem A.4.2 in \cite{bkrw93}.

Define $\Pi_j$ and $\Pi_a$ as the projection operators
\begin{eqnarray*}
\Pi_j g\mapsto\frac{E[g(V,W)Q_{\theta_0}^2|W_j=w_j]}{
E[Q_{\theta_0}^2|W_j=w_j]}, ~~~~
\Pi_a g\mapsto\frac{E[g(V,W)Q_{\theta_0}^2|V=v]}{
E[Q_{\theta_0}^2|V=v]},
\end{eqnarray*}
respectively. Define
\begin{eqnarray*}
D(v,w)&=&\frac{E[ZQ_{\theta_0}^2|V=v,W=w]}{E[Q_{\theta_0}^2|V=v,W=w]},~~~
S(v,w_j)=\frac{E[Q_{\theta_0}^2|V=v,W_j=w_j]}{E[Q_{\theta_0}^2|W_j=w_j]},\\
T(w_i,w_j)&=&\frac{E[Q_{\theta_0}^2|W_i=w_i,W_j=w_j]}{E[Q_{\theta_0}^2|W_j=w_j]},~~~
U(w_j,v)=\frac{E[Q_{\theta_0}^2|W_j=w_j,V=v]}{E[Q_{\theta_0}^2|V=v]}.
\end{eqnarray*}
We say a function $f(s,t)$ belongs to a uniform H\"{o}lder ball $\mathbf{H}_c^r(\mathcal{S}\times\mathcal{T})$ in $t$ relative to $s$ if it is
$J<r$ continuously differentiable w.r.t. $t$ and its $J$-th partial derivative satisfies, with $\kappa\equiv r-J$,
$$\sup_{s\in\mathcal{S}}\sup_{t_1\neq t_2}\frac{|f^{(J)}_t(s,t_1)-f^{(J)}_t(s,t_2)|}{|t_1-t_2|^\kappa}\leq c.$$ Define
$Sf(v,w_j)=S(v,w_j)f_{V|W_j}(v,w_j)$, $Tf(w_i,w_j)=T(w_i,w_j)f_{W_i|W_j}(w_i,w_j)$ and $Uf(w_j,v)=U(w_j,v)f_{W_j|V}(w_j,v)$, where $f_{V|W_j}$,
$f_{W_i|W_j}$ and $f_{W_j|V}$ are the conditional densities of $V$ given $W_j$, $W_i$ given $W_j$ and $W_j$ given $V$ w.r.t. Lebesgue measure,
respectively.

Here, we assume some model assumptions implying that both $b_{jk}^\dag$ and $a_k^\dag$ belong to some H\"{o}lder balls for any $j=1,\ldots,d$ and $k=1,\ldots,l$.
\begin{enumerate}
\item[M5.] We assume that $[\Pi_j D(v,w)]_k\in\mathbf{H}_{\bar c_{j}}^{r_j}[0,1]$,
$Sf(v,w_j)\in\mathbf{H}^{r_j}_{\bar c_j}([l_v,u_v]\times[0,1])$ in $w_j$ relative to $v$ and $Tf(w_i,w_j)\in
\mathbf{H}^{r_j}_{\bar c_{j}}[0,1]^2$ in $w_j$ relative to $w_i$ for  some $0<\bar c_j<\infty$ and $j=1,\ldots,d$.

\item[M6.] We assume that $[\Pi_a D(v,w)]_k\in\mathbf{H}_{\bar c_{0}}^{r_0+1}[l_v,u_v]$ and
$Uf(w_j,v)\in\mathbf{H}^{r_0+1}_{\bar c_0}([0,1]\times[l_v,u_v])$ in $v$ relative to $w_j$ for some $0<\bar c_0<\infty$.
\end{enumerate}
Note that we can simplify $Sf(v,w_j)$ $(Tf(w_i,w_j))$ to $S(v,w_j)$ $(T(w_i,w_j))$ in Condition M5 and simplify $Uf(w_j,v)$ to $U(w_j,v)$ in Condition M6 when we assume that $V$ and $W$ are independent and that $W$ is pairwise independent.

\begin{theorem}\label{asymnor}
Suppose that Conditions M1-M6 and P1-P2 hold. If $K_j\asymp n^{1/(2r_j+1)}$ and $\widetilde I_0$ is invertible, then we have
\begin{eqnarray}
\sqrt{n}(\widehat\beta-\beta_0)=\frac{1}{\sqrt{n}}\sum_{i=1}^{n}\widetilde I_0^{-1}\widetilde\ell_{\beta_0}(X_i)+o_P(1)\overset{d}{\longrightarrow} N(0,\widetilde I_0^{-1}),\label{asynorfor}
\end{eqnarray}
where $\widetilde I_0$ is the efficient information matrix defined as $E\widetilde\ell_{\beta_0}\widetilde\ell_{\beta_0}'$.
\end{theorem}

\section{B-spline Estimate of the Efficient Information}\label{varest}
In this section, we give an explicit B-spline estimate for the efficient information as a by-product of the establishment of asymptotic normality
of $\widehat\beta$. Indeed, it is simply the observed information matrix if we treat the semiparametric model as a parametric one after the
B-spline approximation, i.e., $\mathcal{H}_j=\mathcal{H}_{jn}$ and $\mathcal{G}=\mathcal{G}_n$. Specifically, we treat $\ell_i(\alpha)$ defined in (\ref{sie-lik}) as if it were a parametric likelihood $\ell_i(\beta,\gamma_0, \gamma_1,\ldots,\gamma_d)$.

We construct the corresponding information estimator for $(\beta',\gamma_0,\gamma_1,\ldots,\gamma_2)'$:
\begin{eqnarray*}
\widehat J=\begin{pmatrix}\widehat I_{11} & \widehat I_{12}\\ \widehat I_{21} &\widehat I_{22}\end{pmatrix}_{(l+\sum_{j=0}^d K_j)\times(l+\sum_{j=0}^d K_j)},
\end{eqnarray*}
where $\widehat I_{j,k}=\sum_{i=1}^{n} A_{j}(X_i;\widehat\alpha)A_{k}'(X_i;\widehat\alpha)/n$, for $j,k=1,2$, and
\begin{eqnarray*}
A_1(X;\alpha)&=&\dot\ell_\beta(X;\alpha),\\
A_2(X;\alpha)&=&\left(\dot\ell_{g}[B_{01}],\ldots,\dot\ell_{g}[B_{0K_0}], \dot\ell_{h_1}[B_{11}],\ldots,\dot\ell_{h_d}[B_{dK_d}]\right)'.
\end{eqnarray*}
The parametric inferences imply that the information estimator for $\beta$ is of the form
\begin{eqnarray}
\widehat I=\widehat I_{11}-\widehat I_{12}
\widehat I_{22}^{-1}
\widehat I_{21}.\label{infoest}
\end{eqnarray}
Some calculations further reveal that
\begin{eqnarray}
\widehat I=\mathbb{P}_n\left[\dot\ell_{\widehat\beta}-\dot\ell_{\widehat g}[(\bar\gamma_0^\dag)'\mathbf{B}_0]-
\sum_{j=1}^d\dot\ell_{\widehat h_j}[(\bar\gamma_j^\dag)'\mathbf{B}_j]
\right]^{\otimes 2},\label{ihat2}
\end{eqnarray}
where $[\bar\gamma_j^\dag]_{K_j\times l}=(\gamma_{j1}^\dag,\ldots,\gamma_{jl}^\dag)$ for $j=0,1,\ldots,d$ and
$(\gamma_{0k}^\dag, \ldots, \gamma_{dk}^\dag)^T=\widehat I_{22}^{-1}\widehat I_{21}1_{k}$
where $1_k$ represents the $l$-vector with its $k$-th element as one and others as zeros.
We will use (\ref{infoest}) as our estimator for $\widetilde I_0$.

We need the following additional assumption for Theorem~\ref{infocon}.
\begin{itemize}
\item[M7.] We assume that $$E\sup_{a_k\in\mathcal{G}_n}\left[\int_{l_V}^V[\exp(g(s))-\exp(g_0(s))]a_k(s)ds\right]^2\aplt\|H-H_0\|_2^2.$$
\end{itemize}

\begin{theorem}\label{infocon}
Under Conditions M1-M7 and P1-P2, we have $\widehat I\overset{P}{\rightarrow}\widetilde I_0$.
\end{theorem}

\section{Numerical Results}

\subsection{Simulations}\label{simulat}
We perform a Monte-Carlo study to assess the finite-sample performance of our proposed method.
To compare with the penalized NPMLE in \cite{mk05b}, we adopt the same setting used in their paper.
We simulate the current status data from the partly linear additive Cox model which is a special case of general transformation model.
We choose $H(u) = \log A(e^u)$ where $A(u) = e^{k_0}(\exp(u/3)-1)$ with $k_0=0.06516$. The errors $\epsilon$ follow an extreme
value distribution with $F(s)=1-\exp(-e^s)$. The regression coefficients $\beta_1=0.3$ and $\beta_2 = 0.25$. The covariate $Z_1$ is Uniform$[0.5, 1.5]$ and $Z_2$ is Bernoulli with success probability $0.5$. We choose $W$ as Uniform$[1, 10]$ and $h(w)=\sin(w/1.2-1)-k_0$. Censoring times are standard exponential distribution conditional on being in the interval $[0.2, 1.8]$. The sample sizes are $n=400$ and $n=1600$. We simulate $400$ realizations for both sample sizes.

\begin{figure}
\begin{center}
\resizebox{3in}{3in}{\includegraphics{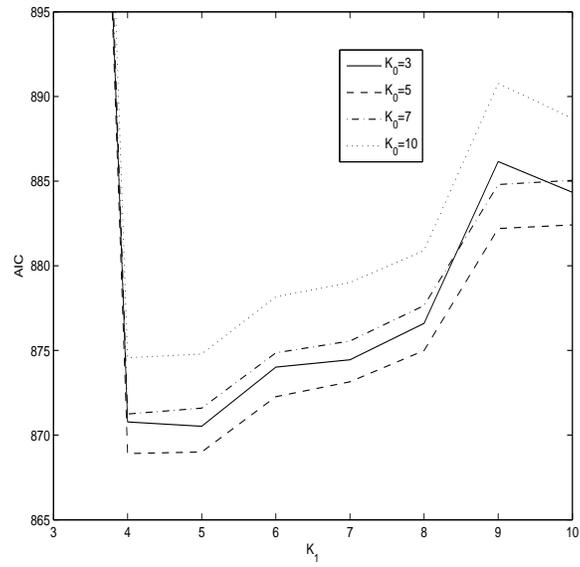}}
\end{center} \caption{AIC scores under different combinations of $K_0$ and $K_1$}\label{fig:aic}
\end{figure}

\begin{figure}
\begin{center}
\resizebox{5in}{2.5in}{\includegraphics{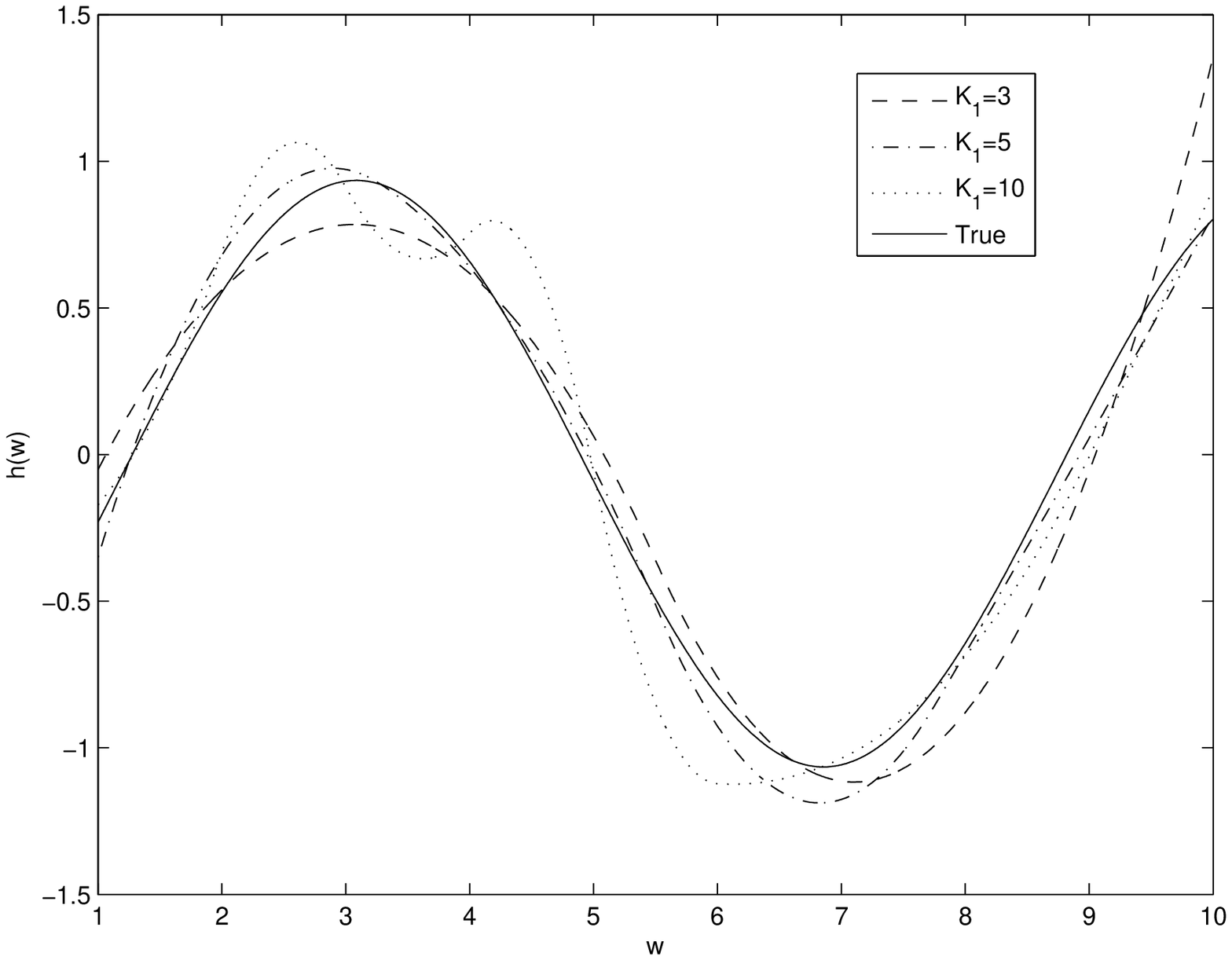}\hspace{.2in}\includegraphics{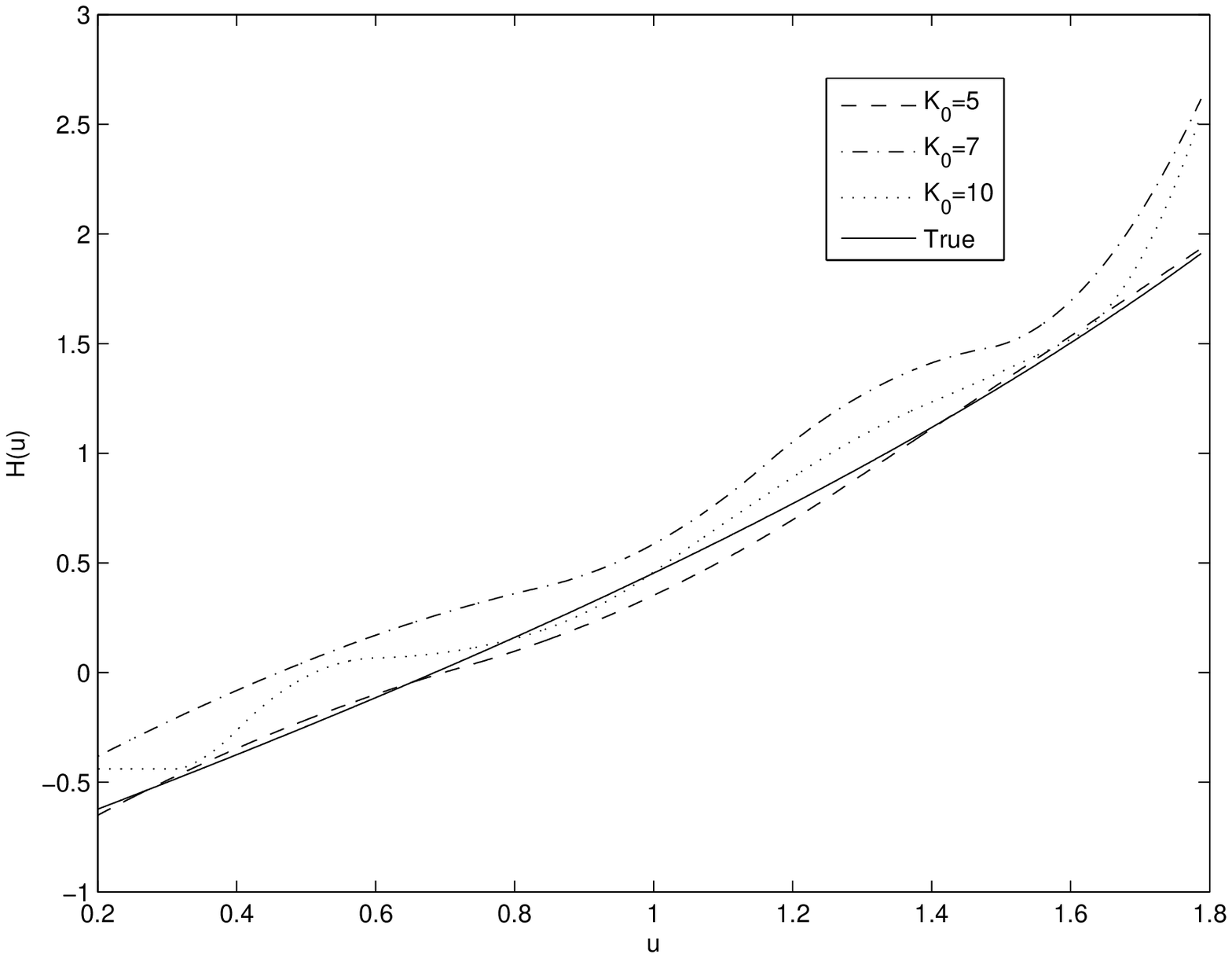}}
\end{center} \caption{Plot of the estimated $h$ and $H$ with various values of $K_0$ and $K_1$.}\label{fig:pic-aic}
\end{figure}

In practice, the numbers of knots for $H$ and $h_j$ need to be determined. Common variable selection methods such as the Akaike information criterion (AIC),
and the Bayesian information criterion (BIC) can be employed for selecting the optimal number of knots. In this paper, we determine $K_0, K_1, \ldots, K_d$ by the AIC given by
$$\mathrm{AIC} = -2 \sum_{i=1}^n \ell_i(\hat\alpha) + 2 (\ell + \sum_{j=0}^d K_j)$$
In our simulation, we use a quadratic spline to approximate both function $h$ and function $g$ in $H$. Then, $\mathrm{AIC}=-2 \sum_{i=1}^n \ell_i(\hat\alpha) + 2(K_0+K_1+2)$. Based on our experiences, it is generally adequate to choose less than ten knots to achieve reasonable approximation, provided that $h$ and $H$ are not overly erratic. Figure \ref{fig:aic} shows the AIC scores under different combinations of $K_0$ and $K_1$ for one realization of the simulation with the sample size $n=1600$. It shows that the optimal choices for $K_0$ and $K_1$ are $5$ and $5$, respectively. The estimated $h$ and $H$ with various values of $K_0$ and $K_1$ are plotted in Figure \ref{fig:pic-aic}.  In the left panel of Figure \ref{fig:pic-aic}, we fix $K_0=5$ and plot the estimated $h$ with $K_1=3, 5, 10$.  When $K_1$ is small (e.g., $K_1=3$), there seems be to a big bias in our estimator. On the other hand, when $K_1$ is large (e.g., $K_1=10$), the estimator displays a wiggly behavior. In the right panel of Figure \ref{fig:pic-aic}, we fix $K_1=5$ and plot the estimated $H$ with $K_0=5, 7, 10$. As the number of knots is increasing, the estimated $H$ shows a similar wiggly shape. Hence, the numbers of knots should be chosen with caution.

\begin{table}
\caption{Monte Carlo results for the partly linear Cox model with current status data based on $400$ replicates}
\begin{center}
\begin{tabular}{llcccc}
\hline
       &     &         Sample size $400$  & Sample size $1600$\\
\hline
$\widehat\beta_1$ &   Bias     &  0.0318              &  0.0100               \\
              &   SD       &  0.2919              &  0.1246               \\
              &  ESD    &  0.3102              &  0.1325              \\
              & Coverage  & 0.9620               & 0.9690 \\
$\widehat\beta_2$ &   Bias     &  0.0168              &  0.0074               \\
              &   SD       &  0.1533              &  0.0797               \\
              &  ESD    &  0.1612              &  0.0803                \\
& Coverage  & 0.9710               & 0.9680 \\
\hline
Joint       &  Coverage &  0.9620   &   0.9550\\
              \hline
\end{tabular}
\end{center}\label{table}
SD: Standard error; ESD: Estimated standard error
\end{table}

\begin{figure}
\begin{center}
\resizebox{5in}{2.5in}{\includegraphics{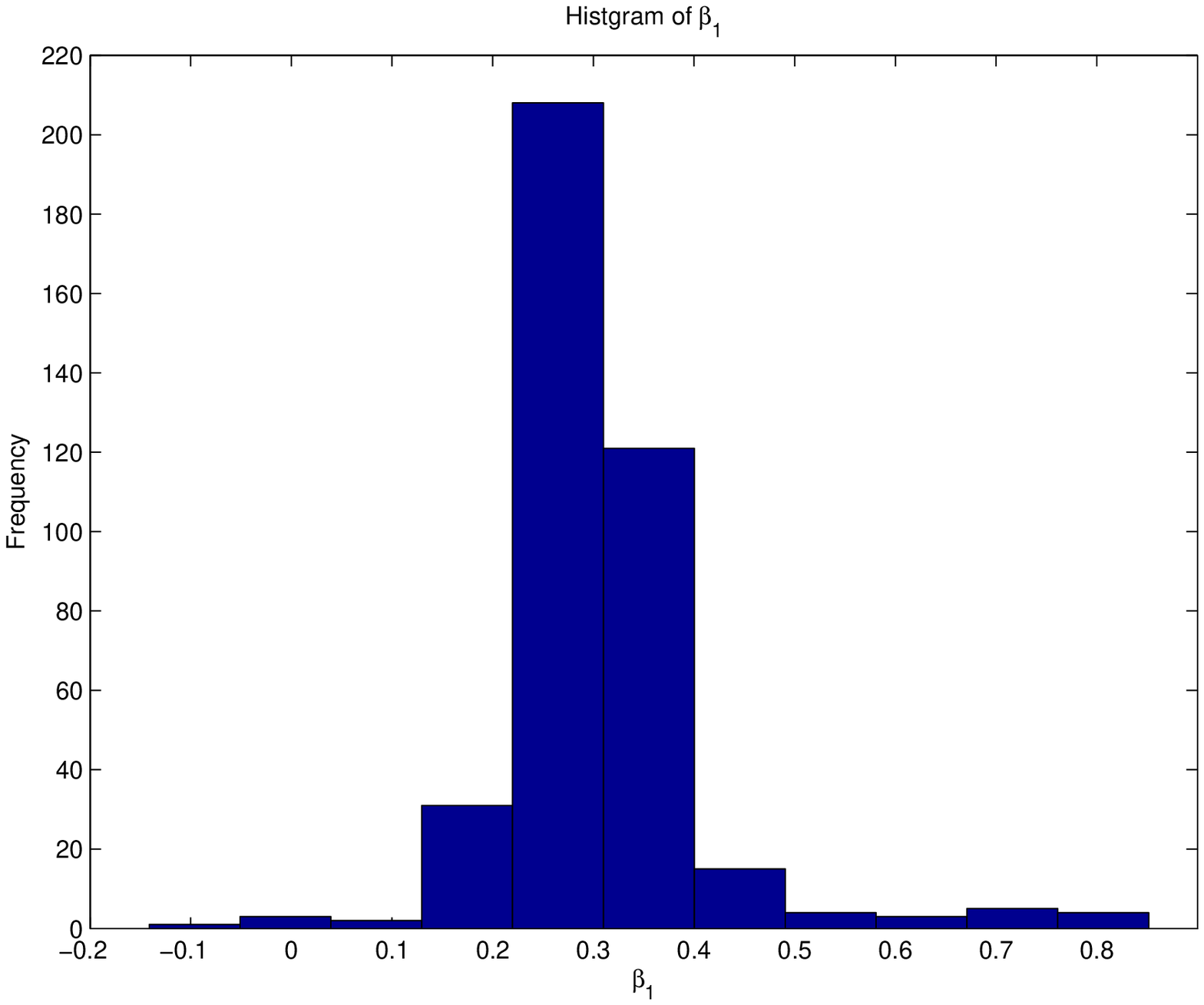}\hspace{.2in}\includegraphics{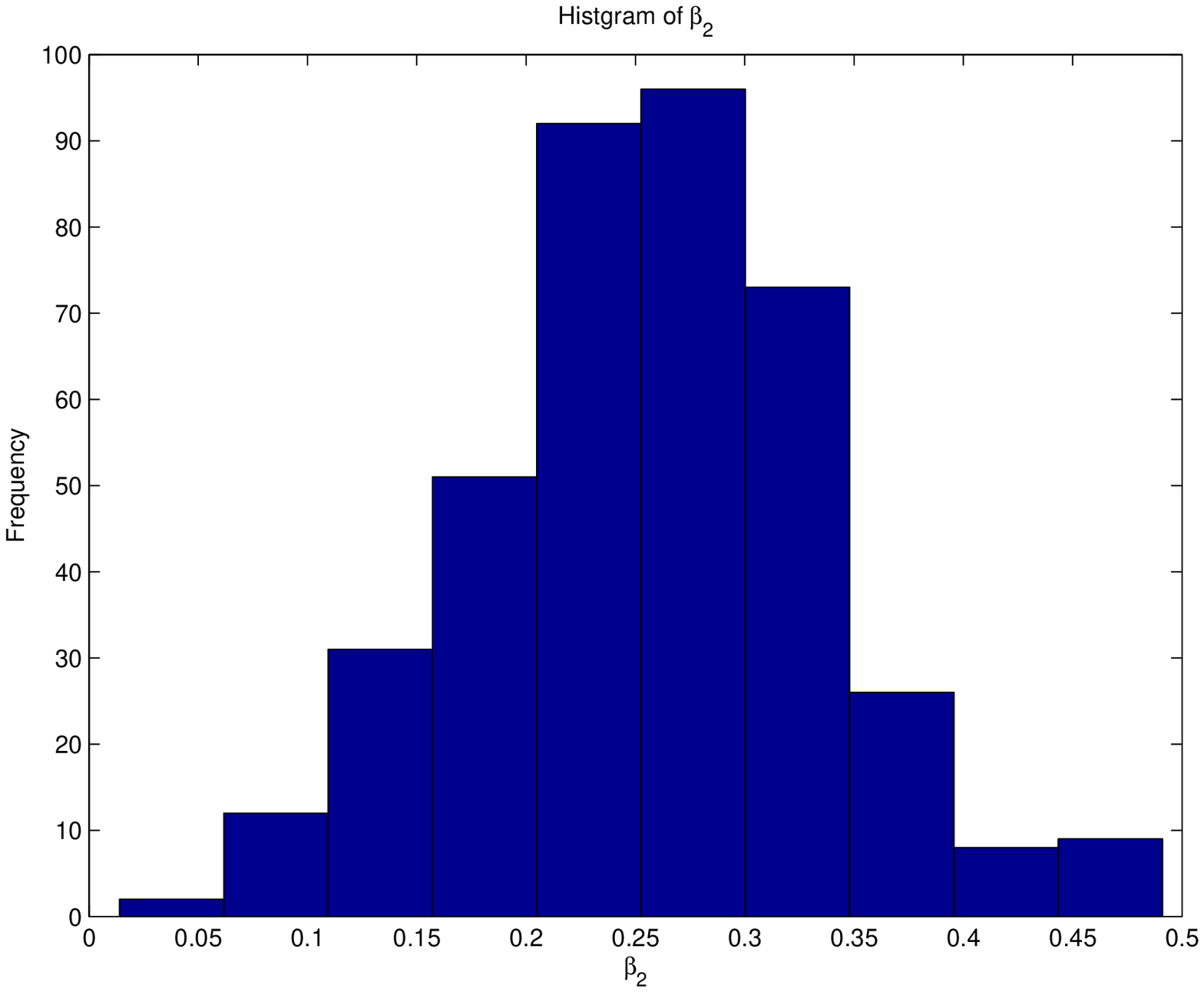}}
\end{center} \caption{Histogram of $\widehat\beta_1$ and $\widehat\beta_2$ based on $1600$ samples and $400$ replicates.}\label{fig:hist}
\end{figure}

\begin{figure}
\begin{center}
\resizebox{5in}{2.5in}{\includegraphics{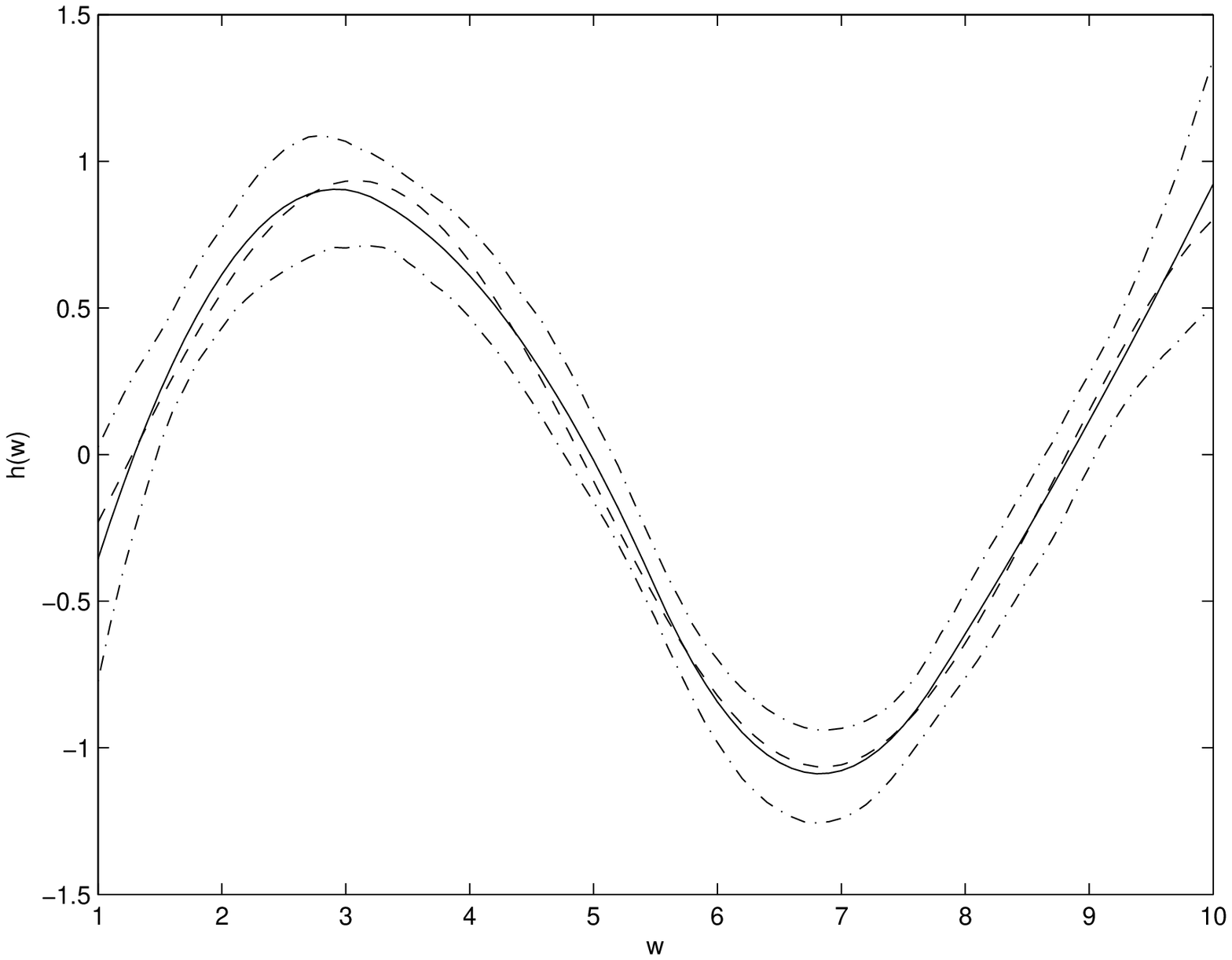}\hspace{.2in}\includegraphics{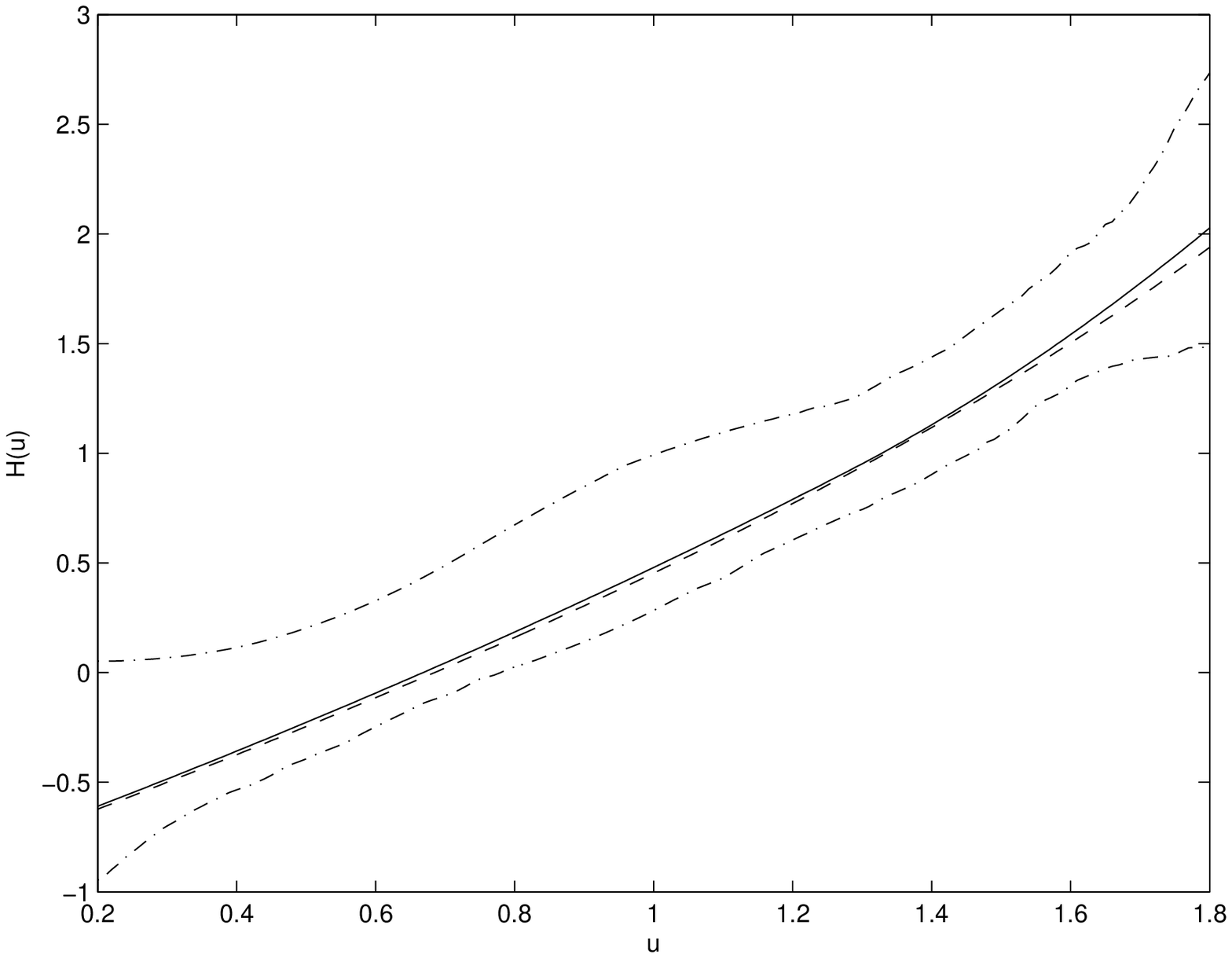}}
\end{center} \caption{Left: Estimate and pointwise confidence interval for $h$. Right: Estimate and pointwise confidence interval for $H$. The solid line is the average
estimate over $400$ realizations from sample size $n=1600$, and the dashed line is the true function. The dash-dotted lines are the $95$\% pointwise confidence interval. }\label{fig:hw}
\end{figure}

Simulation results show that our B-spline estimation procedure performs quite well in the semiparametric transformation model. The bias and standard errors of the spline estimates of $\beta_1$ and $\beta_2$ are given in Table \ref{table}. The table shows that the
sample biases of both $\widehat \beta_1$ and $\widehat\beta_2$ are small. The ratio of the standard errors for the two sample sizes is close to $2$, a result consistent with a $\sqrt{n}$-convergence rate for $\widehat\beta_1$ and $\widehat\beta_2$. The estimated standard errors from (\ref{infoest}) (denoted as ESD) are also displayed in Table \ref{table}, which are very close to the simulation results. Although our proposed method tends to overestimate the standard error slightly but the overestimation lessens as sample size increases. The 95\% confidence interval constructed from (\ref{infoest}) generally have coverage close to the nominal value.
Histograms of $\widehat\beta_1$ and $\widehat\beta_2$ are shown in Figure \ref{fig:hist}. It is clear that the marginal distributions of $\widehat\beta_1$ and $\widehat\beta_2$ are Gaussian.
The left panel of Figure \ref{fig:hw} displays the spline estimate of $h(w)$ and the monotone estimate $\widehat H$ is given
in the right panel of Figure \ref{fig:hw}.
The dashed line is the true function, the solid line is the average estimate over $400$
realizations, and the dash-dotted line is the 95\% pointwise confidence band for $h(w)$ or $H(v)$ when
we know the true model, which is obtained by taking $2.5$ percentile
and $97.5$ percentile of these $400$ estimates at each $w$ or $v$.

To compare our spline based method with the penalized method in \cite{mk05b}, there are four obvious advantages of our method.
First, the computational cost of our spline estimate $\widehat{H}$ is much less expensive than that used in \cite{mk05b}, i.e.
the cumulative sum diagram approach. This is because the number of
basis B-splines (thus the number of knots), e.g.,  $K_0=5$ and $K_1=5$, is often taken much smaller than the
sample size $n$, thus the dimension of the estimation problem is
greatly reduced. Secondly, our estimate of the transformation function $H$ is smooth with a higher convergence rate. We obtain a narrower confidence interval for $H$ shown in the right panel of Figure \ref{fig:hw}. Thirdly, we can obtain an explicit consistent estimate $\widehat I$. However, the block jackknife approach proposed
in \cite{mk05b} is not theoretically justified. At last, we do not require the constrained optimization in our implementations.

\subsection{Application: Calcification data}


\begin{table}
\caption{The estimates and their corresponding estimated standard errors for the parametric part for the calcification data}
\begin{center}
\begin{tabular}{lccccc}
\hline
                &  extreme value distribution & logistic distribution \\ \hline
$\hat\beta_1$   &   $-0.1870$                   &  $ -0.2562$    \\
ESD($\hat\beta_1$) & $0.2322$   & $0.2119$  \\
$\hat\beta_2$   & $0.3502$           & $0.3573$  \\
ESD($\hat\beta_2$) & $0.3481$ &  $0.3280$   \\
\hline
\end{tabular}
\end{center}\label{table2}

ESD: Estimated standard error

\end{table}

We illustrate the proposed method in a dataset from the calcification study. Yu et al. (2001) investigated the calcification of intraocular lenses, which is an infrequently reported complication of cataract treatment. Understanding the effect of some clinical variables on the time to calcification of the lenses after implantation is the objective of the study. The patients were examined  by an ophthalmologist to determine the status of calcification at a random time ranging from zero to thirty six months after implantation of the intraocular lenses. The severity of calcification was graded into five categories ranging from zero to four. In our analysis, we simply treat those with severity $>1$ as calcified and those with severity $\le 1$ as not calcified. This dataset can be treated as the current status dataset because only the examination time and the calcification status at examination are available. The interesting covariates include $Z_1$ incision length, $Z_2$ gender ($0$ for female and $1$ for male), and $W$ age at implantation/10. The original dataset has $379$ records. We remove the one record with missing measurement, resulting the sample size $n=378$. This dataset has been studied by \cite{xue04}, \cite{lam05}, and \cite{ma09}. \cite{xue04} and \cite{lam05} modeled the event time directly and did not use any transformation. A straightforward estimation of the hazard function is not available. \cite{ma09} used the cure model to fit the data, and assumed a generalized linear model for the cure probability. For subjects
not cured, the linear and partly linear Cox proportional hazards models are used
to model the survival risk.

\begin{figure}
\begin{center}
\resizebox{5in}{2.5in}{\includegraphics{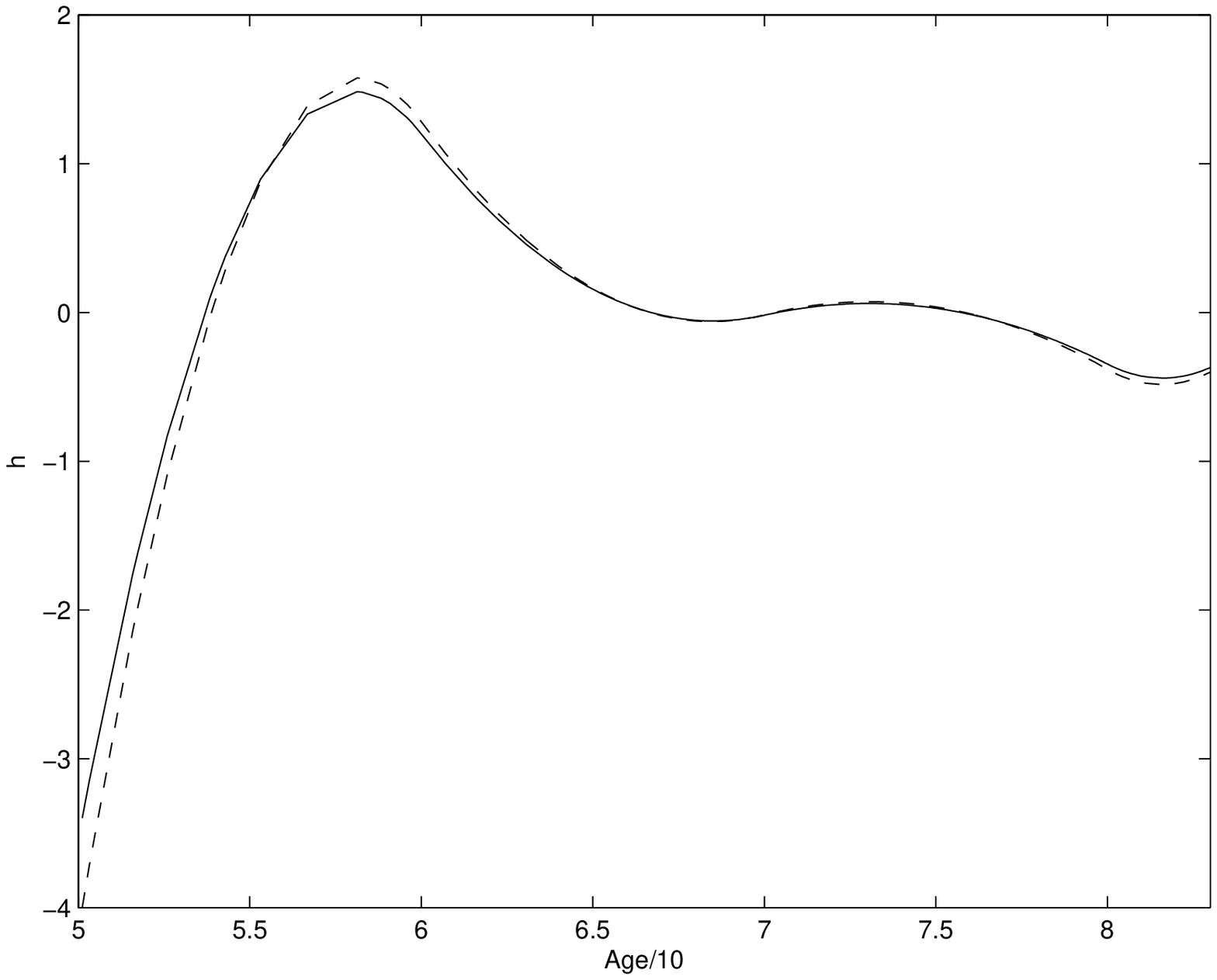}\hspace{.2in}\includegraphics{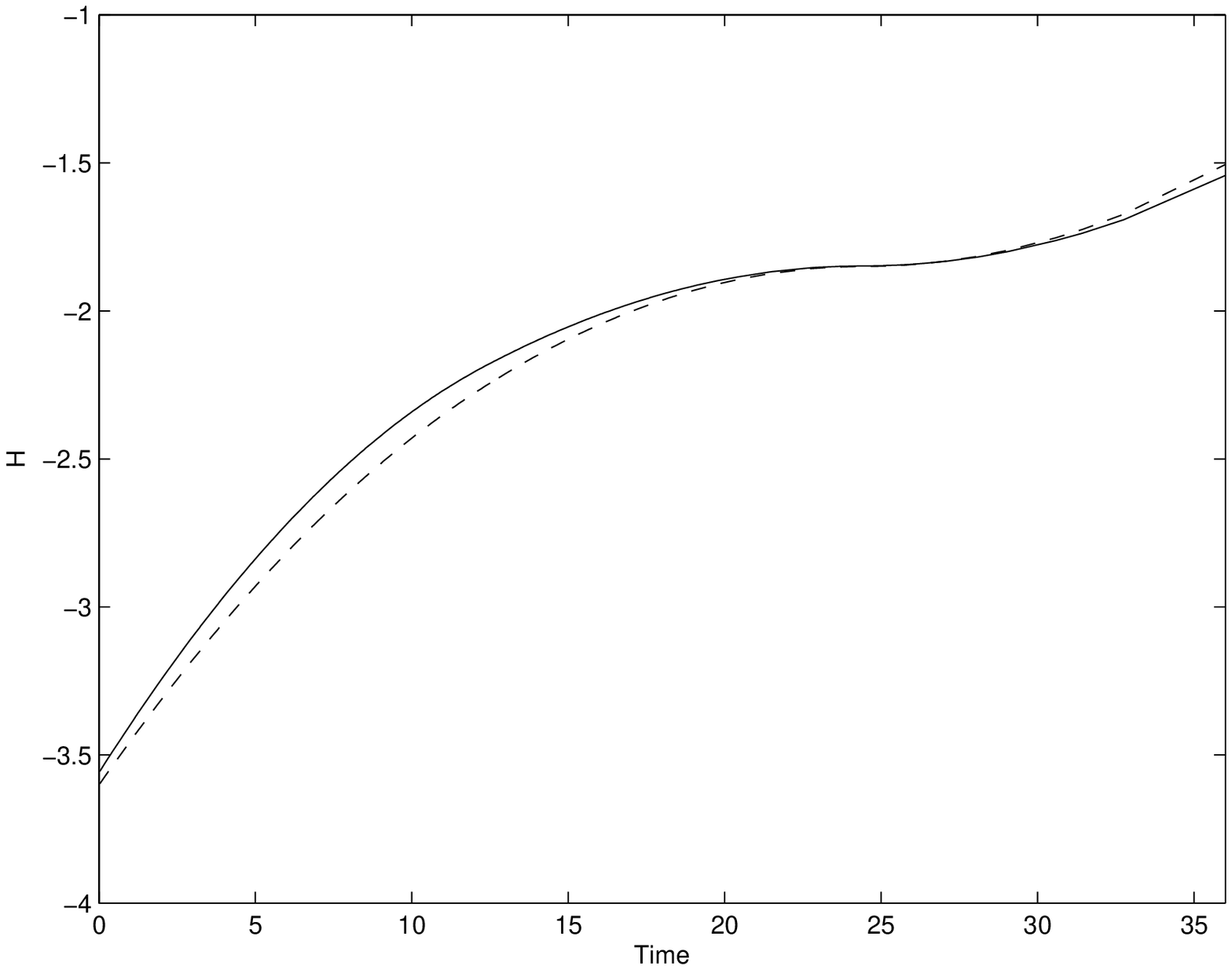}}
\end{center} \caption{The spline estimates of $h(w)$ and $H(v)$ under two different assumptions of the error distribution: extreme value distribution (solid) and logistic distribution (small dashes).}\label{fig:hw-lense}
\end{figure}

We fit this dateset using the semiparametric additive transformation model. We assume the error distribution $F$ to be one of the two distributions: extreme value distribution and logistic distribution. We approximate $h$ and $\log \dot{H}$ by quadratic splines. The optimal choices of knots for $h$ and $\log \dot{H}$ are $6$ and $5$, respectively. The estimates and their corresponding estimated standard errors for the parametric part are summarized in Table \ref{table2}. The estimates for $h(w)$ based on different error distributions are displayed in the left panel of Figure \ref{fig:hw-lense}, and the estimates of $H(v)$ are plotted in the right panel of Figure \ref{fig:hw-lense}. The analysis shows very similar results for these two error distributions.
From Table \ref{table2}, both incision length and gender are insignificant at the 5\% level of significance. From the left panel of Figure \ref{fig:hw-lense}, $h(w)$ increases steadily from age 50, achieving a peak at age 60, decreasing gradually thereafter, which means that patients ages around 60 tend to enjoy a longer time to calcification. The estimated transformation function $\widehat H$ in the right panel of Figure \ref{fig:hw-lense} displays a nonlinear behavior and it shows that the transformation is necessary.

We can incorporate an unknown scale parameter into to the residual error distribution $F(\cdot)$ to further improve the above analysis.
Our general B-spline estimation framework can also handle this type of transformation models easily.

\section*{Acknowledgement}
The first
author's research is supported by the National Science Foundation
under grant DMS-0906497. The second
author's research is supported by the National Science Foundation
under grant CMMI-1030246 and DMS-1042967. The authors would like to thank Professor Alexis K. F. Yu for providing the Calcification data and thank Professors Michael Kosorok and Donglin Zeng for many helpful comments and
suggestions to improve the paper.

\appendix


\setcounter{equation}{0}
\renewcommand{\theequation}{A.\arabic{equation}}

\section*{Appendix}


\subsection*{Some useful Lemmas}
We define $\epsilon$-covering number ($\epsilon$-bracketing number) as $N(\epsilon,\mathcal{A},d)$
($N_{B}(\epsilon,\mathcal{A},d)$). The corresponding $\epsilon$-entropy ($\epsilon$-bracketing entropy)
is defined as $H(\epsilon,\mathcal{A},d)=\log N(\epsilon,\mathcal{A},d)$ ($H_B(\epsilon,\mathcal{A},d)=\log N_{B}(\epsilon,\mathcal{A},d)$).
Define $\mathcal{G}_n(\delta_0; \|\cdot\|)=\{g: g(v)=\gamma_0'\mathbf{B}_0(v)\;\mbox{satisfying}\;\|g\|\leq\delta_0\}$ and $\mathcal{H}_{jn}(\delta_j; \|\cdot\|)=\{h_j: h_j(w_j)=\gamma_j\mathbf{B}_j(w_j)\;\mbox{satisfying}\;\|h_j\|\leq\delta_j\;\mbox{and}\;\int_0^1 h_j(w_j)dw_j=0\}$. Obviously, $\mathcal{G}_n(c_0; \|\cdot\|_\infty)=\mathcal{G}_n$ and $\mathcal{H}_{jn}(c_j;\|\cdot\|_{\infty})=\mathcal{H}_{jn}$.  Lemma~\ref{bsapp} follows from the B-spline approximation property (\ref{sieapp}). Lemma~\ref{entle} is directly implied by Lemma 2.5 in \citep{vg00}. Lemma~\ref{asynorlem} is adapted from Proposition 1 in \citep{ch10}.

\begin{lemma}\label{bsapp}
There exist $g_n\in\mathcal{G}_n$ and $h_{jn}\in\mathcal{H}_{jn}$ such that
\begin{eqnarray}
\|g_n-g_0\|_{\infty}&\asymp&K_0^{-r_0},\label{happrate}\\
\|H_n-H_0\|_{\infty}&=&O(K_0^{-r_0}),\label{bighapprate}\\
\|h_{jn}-h_{j0}\|_{\infty}&\asymp&K_j^{-r_j},\label{gapprate1}\\
\left\|\sum_{j=1}^dh_{jn}-\sum_{j=1}^d h_{j0}\right\|_{\infty}&=&O\left(\max_{j=1,\ldots,d}\{K_j^{-r_j}\}\right),\label{gapprate}
\end{eqnarray}
where $H_n(v)=\int_{l_v}^v\exp(g_n(s))ds$.
\end{lemma}

\begin{lemma}\label{entle}
\begin{eqnarray}
H(\epsilon, \mathcal{G}_n(\delta_0; \|\cdot\|), \|\cdot\|)&\aplt&K_{0}
\log(1+4\delta_0/\epsilon),\label{enth}\\
H(\epsilon, \mathcal{H}_{jn}(\delta_j; \|\cdot\|),\|\cdot\|)&\aplt&K_{j}
\log(1+4\delta_j/\epsilon)\label{ents}
\end{eqnarray}
for $1\leq j\leq d$.
\end{lemma}

\begin{lemma}\label{donle}
Let $\mathbf h=(h_1,\ldots,h_d)$. Define $\mathcal{K}=\{\zeta(\beta, \mathbf h,H): \beta\in\mathcal{B},
\mathbf h\in\prod_{j=1}^d\mathcal{H}_{jn}, g\in\mathcal{G}_n\}$, where the form of $\zeta$
is defined in (\ref{inter4}). We have
\begin{eqnarray}
\sup_{\zeta\in\mathcal{K}}|\mathbb{G}_n\zeta|=O_P(\max_{j=0,1,\ldots,d}\{K_{j}^{1/2}\}).\label{conmod}
\end{eqnarray}
\end{lemma}
{\bf Proof:} Define $l^{\ast}(\beta, \mathbf h, H)=\delta
F(\beta'z+\sum_{j=1}^d h_j(w_j)+H(v))+(1-\delta)[1-F(\beta'z+\sum_{j=1}^d h_j(w_j)+H(v))]$. The construction of $l^\ast(\cdot)$ implies that
\begin{eqnarray}
\|l^{\ast}(\beta_0,\mathbf h_n,H_n)-l^{\ast}(\beta_0,\mathbf h_0,H_0)\|_{\infty}=
O(\max_{j=0,1,\ldots,d}\{K_j^{-r_j}\})\label{lstaapp}
\end{eqnarray}
based on (\ref{bighapprate}), (\ref{gapprate}) and M4. Thus, $l^{\ast}(\beta_0, \mathbf h_n,H_n)$ is bounded away from zero for
sufficiently large $n$.

For any $\beta_1,\beta_2\in\mathcal{B}$,
$\mathbf h_1, \mathbf h_2\in\prod_{j=1}^d\mathcal{H}_{jn}$ and $g_1,g_2\in\mathcal{G}_n$, we have
\begin{eqnarray}
&&|\zeta(\beta_1,\mathbf h_1,H_1)-\zeta(\beta_2, \mathbf h_2,H_2)|\nonumber\\&\aplt&
|l^{\ast}(\beta_1,
\mathbf h_1,H_1)-l^{\ast}(\beta_2, \mathbf h_2,H_2)|\nonumber\\&\aplt&
\|\beta_1-\beta_2\|+\sum_{j=1}^d\|h_{1j}-h_{2j}\|_\infty+\|g_1-g_2\|_\infty.\label{inter5}
\end{eqnarray}
The first and second inequalities in the above follow from the fact that
$l^{\ast}(\beta_0,\mathbf h_n,H_n)$ is strictly positive for sufficiently
large $n$ by (\ref{lstaapp}), and Condition M4(a),
respectively. As shown in (\ref{inter5}), the functions in the class $\mathcal{K}$ are Lipschitz continuous in $(\beta,\mathbf h, g)$.
Therefore, by combining Lemma~\ref{entle} and Theorem 2.7.11 in \citep{vw96}, we obtain that
$$H_{B}(\epsilon,\mathcal{K},L_2(P))\aplt\max_{0\leq j\leq d}\{K_j\}\log(1+M/\epsilon),$$ where $M=\max_{0\leq j\leq d}\{4c_j\}$. In the end, we apply Lemma 3.4.2
in \citep{vw96} to this uniformly bounded class of functions $\mathcal{K}$ to obtain
(\ref{conmod}).  $\Box$

\begin{lemma}\label{asynorlem}
Suppose the following Conditions (B1)-(B3) hold.
\begin{itemize}
\item[B1.] $\mathbb{P}_n\dot\ell_{\widehat\beta}=o_P(n^{-1/2})$, $\mathbb{P}_n\dot\ell_{\widehat g}[\bar a^\dag]=o_P(n^{-1/2})$ and $\mathbb{P}_n\dot\ell_{\widehat h_j}[\bar b_j^\dag]=o_P(n^{-1/2})$;

\item[B2.] $\sup_{\{\alpha: d(\alpha, \alpha_0)\leq C_1n^{-r/(2r+1)}\}}\mathbb{G}_n(\widetilde\ell_{\beta}(X;\alpha)-\widetilde\ell_{\beta}(X;\alpha_0))=o_P(1)$;

\item[B3.] $P(\widetilde\ell_{\beta}(X;\alpha)-\widetilde\ell_{\beta}(X;\alpha_0))=-\widetilde I_0(\beta-\beta_0)+o(\|\beta-\beta_0\|)+
o(n^{-1/2})$ for $\alpha$ satisfying $d(\alpha, \alpha_0)\leq C_1n^{-r/(2r+1)}$.
\end{itemize}
If $\widehat\alpha$ is consistent and $\widetilde I_0$ is invertible, then we have
$$\sqrt{n}(\widehat\beta-\beta_0)=\frac{1}{\sqrt{n}}\sum_{i=1}^{n}\widetilde I_0^{-1}\widetilde\ell_{\beta_0}(X_i)+o_P(1)\overset{d}{\longrightarrow} N(0,\widetilde I_0^{-1}).$$
\end{lemma}

\begin{lemma}\label{holdlem}
(i) If $a(\mathbf{s},t) = a(\mathbf{s}_1,\mathbf{s}_2,t)\in \mathbf{H}_c^r(\mathcal{S}_1\times\mathcal{S}_2\times\mathcal{T})$ in
$t$ relative to $\mathbf{s}_1$ and $\mathbf{s}_2$, then
$\int_{\mathcal{S}_1} a(\mathbf{s}_1,\mathbf{s}_2,t)\, d\mathbf{s}_1\in\mathbf{H}_{c'}^r(\mathcal{S}_2\times\mathcal{T})$ in $t$ relative to $\mathbf{s}_2$.

(ii)
If $a(\mathbf{s},t), b(\mathbf{s},t)\in\mathbf{H}_c^r(\mathcal{S}\times\mathcal{T})$ in $t$ relative to $\mathbf{s}$, then
$c(\mathbf{s},t) \equiv a(\mathbf{s},t)b(\mathbf{s},t)\in\mathbf{H}_{c'}^r(\mathcal{S}\times\mathcal{T})$ in $t$ relative to $\mathbf{s}$.

(iii)
If $a(\mathbf{s},t)\in\mathbf{H}_c^r(\mathcal{S}\times\mathcal{T})$ in $t$
relative to $\mathbf{s}$ and $f(\cdot) \in C^{\lceil \beta
\rceil}$, then $f(a(\mathbf{s},t))\in\mathbf{H}_{c'}^r(\mathcal{S}\times\mathcal{T})$ in $t$ relative to $\mathbf{s}$.
\end{lemma}
{\bf Proof:}
Let $\lfloor r \rfloor$ be the largest integer smaller than $r$. Denote the $m$-th derivative of $a(\mathbf s, t)$ w.r.t. $t$ as $D_t^m a(\mathbf s, t)$ for $m=0,1,\ldots,\lfloor r \rfloor$.

(i) Note that $D_t^{m} a(\mathbf{s}_1,\mathbf{s}_2,t)$ is bounded for
$0\leq m \leq \lfloor r \rfloor$, by the dominated convergence theorem, we can take
derivative inside the integral to obtain
\[
D_t^{m}\Bigl(\int_{\mathcal{S}_1} a(\mathbf{s}_1,\mathbf{s}_2,t)\, d\mathbf{s}_1\Bigr)
=\int_{\mathcal{S}_1} D_t^{m} a(\mathbf{s}_1,\mathbf{s}_2,t)\, d\mathbf{s}_1,
\]
which implies that $D_t^m (\int_{\mathcal{S}_1} a(\mathbf{s}_1,\mathbf{s}_2,t)\,
d\mathbf{s}_1)$ is bounded for $0\leq m\leq \lfloor r \rfloor$. Using this and the
fact that
\begin{align*}
&\frac {\bigl|D_t^{\lfloor r \rfloor}\bigl(\int_{\mathcal{S}_1}
a(\mathbf{s}_1,\mathbf{s}_2,t_2)\, d\mathbf{s}_1\bigr) -
D_t^{\lfloor r \rfloor}\bigl(\int_{\mathcal{S}_1} a(\mathbf{s}_1,\mathbf{s}_2,t_1)\,
d\mathbf{s}_1\bigr)\bigr|}
{|t_2-t_1|^{r-\lfloor r \rfloor}} \\
&\leq \int_{\mathcal{S}_1} \sup_{\mathbf{s}_1, \mathbf{s}_2} \sup_{t_1\ne t_2}
\frac{|D_t^{\lfloor r \rfloor}a(\mathbf{s}_1,\mathbf{s}_2,t_2) -
D_t^{m_{\alpha}}a(\mathbf{s}_1,\mathbf{s}_2, t_1)|}{|t_2 -
t_1|^{r-\lfloor r \rfloor}}\,d\mathbf{s}_1 \leq c'< \infty,
\end{align*}
for all $\mathbf{s}_2$ and $t_1\ne t_2$,
we conclude that $\int_{\mathcal{S}_1} a(\mathbf{s}_1,\mathbf{s}_2,t)\, d\mathbf{s}_1\in\mathbf{H}_{c'}^r(\mathcal{S}_2\times\mathcal{T})$ in $t$ relative to $\mathbf{s}_2$ for some $c'<\infty$.

(ii) The result is true because
\[
D_t^m c = \sum_{i+j=m} D_t^i a D_t^j b
\]
is bounded for $0\leq m \leq \lfloor r\rfloor$. Also we note that for
$i<\lfloor r\rfloor$,
\[
\frac{|D_t^i a(\mathbf{s},t_2) - D_t^i a(\mathbf{s}, t_1)|}{|t_2 -
t_1|^{r-\lfloor r\rfloor}} = \frac{|\int_{t_1}^{t_2} D_t^{i+1} a(\mathbf{s},t)\,
dt|} {|t_2 - t_1|^{r-\lfloor r\rfloor}}.
\]
It can then be easily verified that
\[
\sup_{\mathbf{s}} \sup_{t_1\ne t_2} \frac{|D_t^{\lfloor r \rfloor}c(\mathbf{s},t_2) -
D_t^{\lfloor r\rfloor}c(\mathbf{s}, t_1)|} {|t_2 - t_1|^{r-\lfloor r\rfloor}}< \infty.
\]

(iii) When $0 <\alpha \leq 1$, the result follows from the
observation that
\[
\frac{f(a(\mathbf{s},t_2))-
f(a(\mathbf{s},t_1))}{|t_2-t_1|^\beta} =
\frac{f(a(\mathbf{s},t_2))-
f(a(\mathbf{s},t_1))}{|a(\mathbf{s},t_2)-a(\mathbf{s},t_1)|} \cdot
\frac{|a(\mathbf{s},t_2)-a(\mathbf{s},t_1)|}{|t_2-t_1|^\beta}.
\]
Using the chain rule, the above observation and part (ii) of the
lemma, the desired result can be obtained by induction for general $\beta$. $\Box$

Denote $$S_k(X;\alpha, w_k)=[\dot\ell_{\beta}(X;\alpha)]_k-\dot\ell_g[a_k](X;\alpha)-\sum_{j=1}^d\dot\ell_{h_j}[b_{jk}](X;\alpha),$$
where $w_k=(a_k, b_{1k},\ldots,b_{dk})$. Let $\mathcal{W}_n=\mathcal{G}_n\times\prod_{j=1}^d\mathcal{H}_{jn}$ and $\mathcal{N}_0=\{\alpha\in\mathcal{A}: d(\alpha,\alpha_0)=o(1)\}$.
\begin{lemma}\label{lem4}
Under Conditions M1-M7 \& P1-P2, we have
\begin{eqnarray}
E\sup_{w_k\in\mathcal{W}_n}|S_k(X;\alpha,w_k)-S_k(X;\alpha_0,w_k)|^2\aplt d^2(\alpha, \alpha_0)
\label{m41}
\end{eqnarray}
for all $\alpha\in\mathcal{N}_0$ and $k=1,\ldots,l$.
\end{lemma}
{\bf Proof:} In view of (\ref{scobeta})-(\ref{scohj}) , we can bound the left hand side of (\ref{m41}) by
\begin{eqnarray*}
&\aplt&\|Q_\theta-Q_{\theta_0}\|_2^2+E\left\{\sup_{a_k\in\mathcal{G}_n}\left[\int_{l_v}^V(\exp(g(s))-\exp(g_0(s)))a_k(s)ds\right]^2(Q_\theta-Q_{\theta_0})^2\right\}\\
&&+E\sup_{a_k\in\mathcal{G}_n}\left[\int_{l_v}^V\exp(g_0(s))a_k(s)ds(Q_\theta-Q_{\theta_0})\right]^2\\&&+
E\sup_{a_k\in\mathcal{G}_n}\left[\int_{l_v}^V(\exp(g(s))-\exp(g_0(s)))a_k(s)ds Q_{\theta_0}\right]^2\\
&&+\sum_{j=1}^dE\sup_{b_{jk}\in\mathcal{H}_{jn}}\left[b_{jk}^2(Q_\theta-Q_{\theta_0})^2\right]
\end{eqnarray*}
after some algebra. The compactness of $\mathcal{G}_n$ and $\mathcal{H}_{jn}$ imply that the third and fifth term in the above are both of the order
$\|Q_\theta-Q_{\theta_0}\|_2^2$. For the second term, we can further bound it by
$$E\left[\sup_{a_k\in\mathcal{G}_n}\int_{l_V}^Va_k^2(s)ds\int_{l_V}^V[\exp(g(s))-\exp(g_0(s))]^2ds(Q_{\theta}-Q_{\theta_0})^2\right].$$ Considering the compactness of $\mathcal{G}$ and $\mathcal{G}_n$, we know the second term is also of the order $\|Q_{\theta}-Q_{\theta_0}\|_2^2$. Assumption M4(a) together with Cauchy-Schwartz inequality implies that $\|Q_\theta-Q_{\theta_0}\|_2^2\aplt\|\beta-\beta_0\|^2+\|H-H_0\|_2^2+\|\sum_{j=1}^d (h_j-h_{j0})\|_2^2$. Since we assume that the density for $W$ is bounded away from zero and infinity, we have that $\|\sum_{j=1}^d (h_j-h_{j0})\|_2^2\aplt\sum_{j=1}^d\|h_j-h_{j0}\|_2^2$ considering the identifiability condition $\int_{0}^1 h_j(w_j)dw_j=0$.
Assumption M7 implies that the fourth term is of the order $\|H-H_0\|_2^2$. Considering the form of $d(\alpha, \alpha_0)$, we conclude the whole proof. $\Box$

\subsection*{Proof of Theorem~\ref{consis}}
Recall that $\mathbf h=(h_1,\ldots,h_d)$. Denote $\mathbf h_0$, $\mathbf h_n$ and $\widehat{\mathbf h}$ as the corresponding true value, B-spline approximation and sieve estimate, respectively. Recall that $l^{\ast}(\beta_0, \mathbf h_n,H_n)$ is bounded away from zero for sufficiently large $n$ as implied by (\ref{lstaapp}). Then, by the definition of $\widehat\alpha$, we have
$$\mathbb{P}_n\log\{l^{\ast}(\widehat{\beta},\widehat{\mathbf h},\widehat{H})/l^{\ast}(\beta_0, \mathbf h_n,H_n)\}\geq
0,$$ which implies that, by the inequality that $\alpha\log(x)\leq\log(1+\alpha(x-1))$ for any
$x>0$ and $\alpha\in(0,1)$,
\begin{eqnarray}
0\leq\mathbb{P}_n\log\left[1+\alpha\left\{\frac{l^{\ast}(\widehat{\beta},
\widehat{\mathbf h},\widehat{H})}{l^{\ast}(\beta_0,\mathbf h_n,H_n)}-1
\right\}\right]\equiv\mathbb{P}_n\zeta(\widehat{\beta},\widehat{\mathbf h},
\widehat{H}).\label{inter4}
\end{eqnarray}
Lemma~\ref{donle} implies that
$(\mathbb{P}_n-P)\zeta(\widehat{\beta},\widehat{\mathbf h},\widehat{H})=
o_P(1)$ since $K_j/n=o(1)$ for any $j=0,1,\ldots,d$. Thus,
$P\zeta(\widehat{\beta},\widehat{\mathbf h},\widehat{H})\geq o_P(1)$ based
on (\ref{inter4}). Let
$U_n(X)=l^{\ast}(\widehat{\beta},\widehat{\mathbf h},\widehat{H})/
l^{\ast}(\beta_0,\mathbf h_n,H_n)$. Based on (\ref{lstaapp}) we know
$PU_n(X)=1+o_P(1)$, which further implies
$P\zeta(\widehat{\beta},\widehat{\mathbf h},\widehat{H})\leq o_P(1)$ by the
concavity of $s\mapsto\log(s)$. This in turn implies that
$P\zeta(\widehat{\beta},\widehat{\mathbf h},\widehat{H})=o_P(1)$. This forces $P|(\beta_0'Z+\sum_{j=1}^d h_{jn}(W_j)+H_n(V))-(\widehat\beta'Z+\sum_{j=1}^d \widehat h_{j}(W_j)+\widehat H(V))|=o_P(1)$ by the strict concavity of $s\mapsto \log s$, Conditions M4(a), P1 and P2. It is easy to verify that $ER_n^2=o_P(1)$ if $E|R_n|=o_P(1)$. Thus, we further have
\begin{eqnarray*}
P\left\{(\widehat{\beta}-\beta_0)'Z+\sum_{j=1}^d (\widehat{h}_j-h_{jn})(W_j)+\widehat{H}(V)-
H_n(V)\right\}^2=o_P(1).
\end{eqnarray*}
Combining the above equation with the identifiability condition M3, we can show
$(\widehat{\beta}-\beta_0)=o_P(1)$. This, in turn, implies that
$$P\left\{\sum_{j=1}^d (\widehat{h}_j-h_{jn})(W_j)+\widehat{H}(V)-
H_n(V)\right\}^2=o_P(1).$$ Since we assume that the joint density of $(V,W)$ is bounded away from zero in M2(b), we have
$$\int_{0}^1\cdots\int_{0}^1\int_{l_v}^{u_v}\left\{\sum_{j=1}^d(\widehat h_j-h_{jn})(w_j)+\widehat H(v)-H_n(v)\right\}^2dvdw_1\cdots dw_d=o_P(1).$$ Considering that $\int_{0}^1h_j(w_j)dw_j=0$ for $h_j\in\mathcal{H}_j\cup\mathcal{H}_{jn}$ and that the joint density of $(V,W)$ is bounded away from infinity, we have $\sum_{j=1}^d\|\widehat h_j-h_{jn}\|_2+\|\widehat H-H_n\|_2=o_P(1)$. The spline approximation result (\ref{bighapprate}) and (\ref{gapprate1}) conclude the proof of (\ref{consfor}).

We next verify the conditions of Theorem 3.2.5 in \cite{vw96} to establish the convergence rate result (\ref{convfor}). Recall that $\theta(z,v,w)=\beta'z+H(v)+\sum_{j=1}^d h_j(w_j)$. Denote $\widehat\theta=z'\widehat\beta+\widehat H(v)+\sum_{j=1}^d \widehat h_j(w_j)$ as its sieve estimate. Following similar arguments in proving the consistency, it suffices to show that
\begin{eqnarray}
\|\widehat\theta-\theta_0\|_2=O_P(n^{-r/(2r+1)}),\label{inter00}
\end{eqnarray}
where $r=\min_{0\leq j\leq d}\{r_j\}$. We first need to  show that
\begin{eqnarray}
P[\ell(\alpha_0)-\ell(\alpha)]\apgt\|\theta-\theta_0\|_2^2 \label{con1}
\end{eqnarray}
for every $\alpha$ in the neighborhood of $\alpha_0$. Define $q(\delta,t)=\delta\log(F(t))+(1-\delta)\log (1-F(t))$ and $\ddot q(\delta, t)$ as its second derivative w.r.t. $t$. Since $\alpha_0$ maximizes $\alpha\mapsto P\ell(\alpha)$, we have
$$P[\ell(\alpha_0)-\ell(\alpha)]=P\left[\frac{-\ddot q(\delta,\widetilde\theta)}{2}(\theta-\theta_0)^2\right],$$ where
$\widetilde\theta$ is on the line segment between $\theta$ and $\theta_0$. The compactness of the parameter spaces imply
that $P[\ell(\alpha_0)-\ell(\alpha)]\asymp\|\theta-\theta_0\|_2^2$. This completes the proof of (\ref{con1}). We next calculate
the order of $E\sup_{\|\theta-\theta_0\|_2\leq\delta}|\mathbb{G}_n(\ell(\alpha)-\ell(\alpha_0))|$ as a function of $\delta$, denoted
as $\phi_n(\delta)$, by the use of Lemma 3.4.2 of \cite{vw96}. Let $\mathcal{F}_{1n}(\delta)=\{\ell(\alpha)-\ell(\alpha_0):
g\in\mathcal{G}_n, h_{j}\in\mathcal{H}_{jn}, \|\theta-\theta_0\|_2\leq\delta\}$. Using the same argument as that in the proof of Lemma~\ref{donle},
we obtain that $H_B(\epsilon, \mathcal{F}_{1n}(\delta), L_2(P))$ is bounded by $C\max_{0\leq j\leq d}\{K_j\}\log(1+\delta/\epsilon)$. This leads to
$$J_{B}(\delta, \mathcal{F}_{1n}(\delta), L_2(P))=\int_{0}^\delta\sqrt{1+H_B(\epsilon, \mathcal{F}_{1n}(\delta), L_2(P))}d\epsilon\leq C\max_{0
\leq j\leq d}\{\sqrt{K_j}\}\delta.$$ The compactness of $\mathcal{G}_n$ and $\mathcal{H}_{jn}$ implies the uniform boundedness of any
$f\in\mathcal{F}_{1n}(\delta)$. Thus, Lemma 3.4.2 of \cite{vw96} gives $\phi_n(\delta)=\max_{0\leq j\leq d}\{\sqrt{K_j}\}\delta+\max_{0
\leq j\leq d}\{K_j\}/\sqrt{n}$. By solving $\delta_{1n}^{-2}\phi_n(\delta_{1n})\leq\sqrt{n}$, we get \begin{eqnarray}
\delta_{1n}=O(\max_{0\leq j\leq d}
\{\sqrt{K_j}\}/\sqrt{n}).\label{delta1n}
\end{eqnarray}

In the end, we show that $\mathbb{P}_n\ell(\widehat\alpha)-\mathbb{P}_n\ell(\alpha_0)\geq -O_P(\delta_{2n}^2)$, where $\delta_{2n}=\max_{0\leq j\leq d}\{K_j^{-r_j}\}$. The definition of $\widehat\alpha$ implies that
$$\mathbb{P}_n[\ell(\widehat\alpha)-\ell(\alpha_0)]\geq A_n+B_n,$$ where $A_n=(\mathbb{P}_n-P)\{\ell(\beta_0,H_n,\mathbf h_n)-\ell(\alpha_0)\}$ and $B_n=P\{\ell(\beta_0,H_n,\mathbf h_n)-\ell(\alpha_0)\}$. A straightforward Taylor expansion gives
\begin{eqnarray*}
A_n=(\mathbb{P}_n-P)\left\{\dot\ell_2(\beta_0, \widetilde H_n, \widetilde{\mathbf h}_n)(H_n-H_0)+\sum_{j=1}^d\dot\ell_{2+j}(\beta_0, \widetilde H_n, \widetilde{\mathbf h}_n)(h_{jn}-h_{j0})\right\},
\end{eqnarray*}
where $\dot\ell_{t}$ is the Fr\'{e}chet derivative of $\ell(\beta_0, H_n, \mathbf h_n)$ w.r.t. the $t$-th argument. Considering (\ref{bighapprate}), (\ref{gapprate1}) and the fact that $0<\epsilon_1\leq|\dot q(\delta, t)|\leq \epsilon_2<\infty$ for $t$ in some compacta of $\mathbb{R}^1$, we have
\begin{equation}\label{inter01}
\begin{split}
P\left\{\frac{\dot\ell_2(\beta_0, \widetilde H_n, \widetilde{\mathbf h}_n)(H_n-H_0)+\sum_{j=1}^d\dot\ell_{2+j}(\beta_0, \widetilde H_n, \widetilde{\mathbf h}_n)(h_{jn}-h_{j0})}{\max_{0\leq j\leq d}\{K_j^{-r_j}\}n^\epsilon}\right\}^2\rightarrow 0
\end{split}
\end{equation}
for any $\epsilon>0$. Let $\mathcal{F}_{2n}=\{\ell(\beta_0, H, \mathbf h)-\ell(\alpha_0): g\in\mathcal{G}_n, h_{j}\in\mathcal{H}_{jn}, \|g-g_0\|_\infty\leq C_0 K_0^{-r_0}, \|h_j-h_{j0}\|_\infty\leq C_jK_j^{-r_j}\}$. Similar analysis in Lemma~\ref{donle}
show that the bracketing entropy integral (in terms of $L_2(P)$) for $\mathcal{F}_{2n}$ is finite, thus yields that $\mathcal{F}_{2n}$ is P-Donsker. Combining this P-Donsker result and (\ref{inter01}), we use Corollary 2.3.12 of \cite{vw96} to conclude that $\sqrt{n}A_n/(\max_{0\leq j\leq d}\{K_j^{-r_j}\}n^\epsilon)=o_P(1)$. By choosing some proper $0<\epsilon<1/2$ satisfying $n^{\epsilon-1/2}=\max_{0\leq j\leq d}\{K_j^{-r_j}\}$, we have $A_n=o_P(\max_{0\leq j\leq d}\{K_j^{-2r_j}\})$. We can also show $B_n\geq -O(\max_{0\leq j\leq d}\{K_j^{-2r_j}\})$ by similar analysis of (\ref{con1}). This shows that
\begin{eqnarray}
\delta_{2n}=\max_{0\leq j\leq d}\{K_j^{-r_j}\}.\label{delta2n}
\end{eqnarray}

Therefore, we have that $d(\widehat\alpha, \alpha_0)=O_P(\delta_{1n}\vee\delta_{2n})$, i.e., (\ref{convfor0}), which directly implies (\ref{convfor})
by choosing $K_j\asymp n^{1/(2r_j+1)}$. $\Box$

\subsection*{Proof of Theorem~\ref{asymnor}}
We apply Lemma~\ref{asynorlem} to prove this theorem. We first check Condition B1. Obviously, $\mathbb{P}_n\dot\ell_{\widehat\beta}=0$ since $\widehat\beta$ maximizes $l(\beta, \widehat g,\widehat h_1,\ldots,\widehat h_d)$, $\widehat\beta$ is consistent and $\beta_0$ is an interior point of $\mathcal{B}$.
Following the analysis in Page
2282 of \cite{mk05b}, we can write, with $\bar a^\dag_I(v)=\int_{l_v}^v\exp(g_0(s))\bar a^\dag(s)ds$,
\begin{eqnarray*}
\bar b_j^\dag&=&\Pi_j D(v,w)-\Pi_j\bar a_I^\dag(v)-\sum_{i\neq j}\Pi_j \bar b_i^\dag\\
&=&\Pi_j D(v,w)-\int_{l_v}^{u_v}\bar a_I^\dag(v)Sf(v,w_j)dv-\sum_{i\neq j}\int_{0}^{1}\bar b_i^\dag(w_i)Tf(w_i,w_j)dw_i.
\end{eqnarray*}
According to Lemma~\ref{holdlem} and dominated convergence theorem, we know that $b_{jk}^\dag(w_j)\in\mathbf{H}^{r_j}_{\widetilde c_j}[0,1]$ under Condition M5,
$b_{jk}^\dag\in L_2^0(w_j)$ and $a^\dag_k\in L_2(H)$ (thus $a_{Ik}^\dag$ is uniformly bounded) for some $0<\widetilde c_j<\infty$.  Then, for each $b_{jk}^\dag$, there exists a $b^\dag_{jkn}\in\mathcal{H}_{jn}$ such that
\begin{eqnarray}
\|b_{jk}^\dag-
b_{jkn}^\dag\|_\infty=O(n^{-r_j/(2r_j+1)})\label{lfsapp}
\end{eqnarray}
by (\ref{sieapp}) and the assumption that $K_j\asymp n^{1/(2r_j+1)}$.

Since $\mathbb{P}_n\dot\ell_{\widehat h_j}[b_{jkn}]=0$ for any $b_{jkn}\in\mathcal{H}_{jn}$,
it suffices to show that
\begin{eqnarray}
\mathbb{P}_n\left\{\dot\ell_{\widehat h_j}[b_{jkn}^\dag]-\dot\ell_{\widehat h_j}[b_{jk}^\dag]\right\}=o_P(n^{-1/2}).\label{inter1}
\end{eqnarray}
We can decompose the left hand side of (\ref{inter1}) as $I_{1n}+I_{2n}$, where
\begin{eqnarray*}
I_{1n}&=&P\left\{\dot\ell_{\widehat h_j}[b_{jkn}^\dag-b_{jk}^\dag]-\dot\ell_{h_{j0}}[b_{jkn}^\dag-b_{jk}^\dag]\right\},\\
I_{2n}&=&(\mathbb{P}_n-P)\left\{\dot\ell_{\widehat h_j}[b_{jkn}^\dag-b_{jk}^\dag]\right\}.
\end{eqnarray*}
By Cauchy-Schwartz Inequality, we have $I_{1n}\aplt\|b_{kjn}^\dag-b_{kj}^\dag\|_\infty\|\widehat\theta-\theta_0\|_2$ based on Conditions M4(a), P1 \& P2. Thus, (\ref{inter00}) and (\ref{lfsapp}) imply that $I_{1n}=O_P(n^{-2r/(2r+1)})=o_P(n^{-1/2})$ since $r>1/2$.
Define $\mathcal{A}_n(\delta)=\{\alpha\in\mathcal{A}_n: d(\alpha, \alpha_0)\leq C_1\delta\}$ and $\mathcal{H}_{jn}'(\delta)=\{b_{jkn}\in\mathcal{H}_{jn}:
\|b_{jkn}-b_{jk}^\dag\|_\infty\leq C_2\delta\}$ for some $0<C_1,C_2<\infty$. As for
the term $I_{2n}$, we first consider the following class of functions:
$$\mathcal{I}_n=\left\{\dot\ell_{h_j}[b_{jkn}-b_{jk}^\dag](X;\alpha): \alpha\in\mathcal{A}_n(n^{\frac{-r}{2r+1}})\;\mbox{and}
\;b_{jkn}\in\mathcal{H}_{jn}'(n^{\frac{-r_j}{2r_j+1}})\right\}.$$ For simplicity,
we write the function in $\mathcal{I}_n$ as $f_{\theta, b_{jkn}}(x)$. Let $\Theta_n(\delta)=\{\beta'z+H(v)+\sum_{j=1}^dh_j(w_j): \alpha\in\mathcal{A}_n(\delta)\}$. It is easy to verify that, for every $x$,
\begin{eqnarray}
\;\;\;\;\;\;\;\;\;\;\;\;|f_{\theta_1,b_{jkn1}}(x)-f_{\theta_2,b_{jkn2}}(x)|\aplt \|\theta_1-\theta_2\|_\infty+\|b_{jkn1}-b_{jkn2}\|_\infty,\label{inter3}
\end{eqnarray}
where $\theta_j\in\Theta_n(n^{-r/(2r+1)})$ for $j=1,2$.
Let $\theta^1,\ldots,\theta^{N(\epsilon, \Theta_n(n^{-r/(2r+1)}), \|\cdot\|_\infty)}$ and $b_{jkn}^1,\ldots,b_{jkn}^{N(\epsilon, \mathcal{H}_{jn}'
(n^{-r_j/(2r_j+1)}), \|\cdot\|_\infty)}$ be
the $\epsilon$-cover for $\Theta_n(n^{-r/(2r+1)})$ and $\mathcal{H}_{jn}'(n^{-r_j/(2r_j+1)})$, respectively. Thus, we can construct the bracket $[f_{\theta^i,b_{jkn}^l}-2C\epsilon,
f_{\theta^i,b_{jkn}^l}+2C\epsilon]$ covering $\mathcal{I}_n$. The bracket size is $4C\epsilon$. Hence, we obtain
\begin{eqnarray*}
&&H_B(\epsilon, \mathcal{I}_n, L_2(P_X))\\
&\leq&H(\epsilon/(4C),\Theta_n(n^{\frac{-r}{2r+1}}), \|\cdot\|_\infty)+
H(\epsilon/(4C),\mathcal{H}_{jn}'(n^{\frac{-r_j}{2r_j+1}}), \|\cdot\|_\infty)\\
&\aplt&\max_{0\leq j\leq d}\{K_j\}\log(1+n^{-r/(2r+1)}/\epsilon)
\end{eqnarray*}
based on Lemma~\ref{entle}. We next apply Lemma~3.4.2 in \cite{vw96} to show $E\|\mathbb{G}_n\|_{\mathcal{I}_n}=o(1)$ which yields $I_{2n}=o_P(n^{-1/2})$.
We first calculate the $\delta$-bracketing entropy integral
$$J_B(\delta,\mathcal{I}_n, L_2(P_X))\equiv\int_0^\delta\sqrt{1+H_B(\epsilon, \mathcal{I}_n,L_2(P_X))}=\max_{0\leq j\leq d}\{\sqrt{K_j}\}
n^{-\frac{r}{4r+2}}\delta^{1/2}.$$ Note that $\|f\|_2\aplt\|b_{jkn}-b{jk}^\dag\|_2$ and $\|f\|_\infty\leq\|b_{jkn}-b{jk}^\dag\|_\infty$ for any $f\in\mathcal{I}_n$, and thus
$\delta$ and $M$ in Lemma 3.4.2 of \cite{vw96} are both chosen as $K_j^{-r_j}$, i.e., $n^{-r_j/(2r_j+1)}$. Then, by Lemma 3.4.2 of \cite{vw96} and some algebra, we have that $$E\|\mathbb{G}_n\|_{\mathcal{I}_n}=O\left(n^{-\left(\frac{r-1}{4r+2}+\frac{r_j}{4r_j+2}\right)}\vee n^{-\frac{4r-1}{4r+2}}\right)=o(1).$$
We have thus verified that $\mathbb{P}_n\dot\ell_{\widehat h_j}(\bar b_j^\dag)=o_P(n^{-1/2})$.

We next show that $\mathbb{P}_n\dot\ell_{\widehat g}[\bar a^\dag]=o_P(n^{-1/2})$ by similar arguments. Similarly, we have
$$\bar a^\dag_I(v)=\Pi_a D(v,w)-\sum_{j=1}^d\int_{0}^1\bar b_j^\dag(w_j)Uf(w_j,v)dw_j.$$ Recall that $\bar a^\dag_I(v)=\int_{l_v}^v\exp(g_0(s))\bar a^\dag(s)ds$. Under Condition M6 and the assumption that $g_0\in\mathbf{H}^{r_0}_{c_0}
[l_v,u_v]$, we can show
that $a_{Ik}^\dag\in\mathbf{H}^{r_0+1}_{\widetilde c_0}[l_v,u_v]$, which implies that $a_k^\dag\in\mathbf{H}^{r_0}_{\widetilde c_0}[l_v,u_v]$ for some
$0<\widetilde c_0<\infty$, based on Lemma~\ref{holdlem}. We next show that $I_{1n}'=o_P(n^{-1/2})$ and $I_{2n}'=o_P(n^{-1/2})$, where
\begin{eqnarray*}
I_{1n}'&=&P\left\{\dot\ell_{\widehat g}[a_{kn}^\dag-a_{k}^\dag]-\dot\ell_{g_{0}}[a_{kn}^\dag-a_{k}^\dag]\right\},\\
I_{2n}'&=&(\mathbb{P}_n-P)\left\{\dot\ell_{\widehat g}[a_{kn}^\dag-a_{k}^\dag]\right\},
\end{eqnarray*}
and $a_{kn}^\dag\in\mathcal{G}_n$ satisfies $\|a_{kn}^\dag-a_k^\dag\|_\infty=O(K_0^{-r_0})$ for any $k=1,\ldots,l$. Similarly, by Cauchy-Schwartz Inequality, we can show that
\begin{eqnarray*}
I_{1n}'&\aplt&\|a_{kn}^\dag-a_k^\dag\|_\infty\|\widehat\theta-\theta_0\|_2+P\left[\int_{l_v}^v(\exp(\widehat g)-\exp(g_0))(s)(a_{kn}^\dag-a_k^\dag)(s)ds\right]\\
&\aplt&\|a_{kn}^\dag-a_k^\dag\|_\infty\left(\|\widehat\theta-\theta_0\|_2+\|\widehat H-H_0\|_2\right)\\
&\aplt&O_P(n^{-r/(2r+1)})=o_P(n^{-1/2})
\end{eqnarray*}
by choosing $K_j\asymp n^{1/(2r_j+1)}$. Following similar arguments in
analyzing $I_{2n}$, we can show that $I_{2n}'=o_P(n^{-1/2})$. Thus, we have verified Condition B1 in Lemma~\ref{asynorlem}. We again apply Lemma 3.4.2 of \cite{vw96} to verify Assumption B2. The details are skipped due to the similarity of the previous analysis.

It remains to verify Assumption B3. This can be easily established using the Taylor expansion in Banach space. However, we first need to reparameterize the efficient score function $\widetilde\ell_\beta(X;\alpha)$ as
\begin{eqnarray*}
\widetilde\ell_\beta(X;\alpha^\ast)&=&ZQ_{\theta}(X)-\left[\int_{l_v}^V\bar a^\dag(s) dH(s)+\sum_{j=1}^d\bar b_j^\dag(W_j)\right]Q_{\theta}(X)\\
&\equiv&\dot\ell_\beta(X;\alpha^\ast)-\dot\ell_\eta[\bar c^\dag](X;\alpha^\ast),
\end{eqnarray*}
where
$\alpha^\ast=(\beta, H, h_1,\ldots, h_d)$, $\eta=(H, h_1,\ldots, h_d)$ and $\bar c^\dag=(\bar a^\dag, \bar b_1^\dag, \ldots, \bar b_d^\dag)$. We first derive two useful equalities (\ref{effinf})-(\ref{inter06}). Let
$E_{\alpha^\ast}$ be the expectation corresponding to the reparametrized likelihood under the parameter $\alpha^\ast$. Since $E_{\alpha^\ast}\widetilde\ell_{\beta}(X;\alpha^\ast)=0$, we have
\begin{eqnarray}
\frac{\partial}{\partial t}|_{t=0} E_{\alpha_t^\ast}\widetilde\ell_{\beta}(X;\alpha_t^\ast)=0,\label{inter03}
\end{eqnarray}
where $\alpha_t^\ast=\alpha_0^\ast+t\epsilon$. Define $\widetilde\ell_{\beta,\beta}$ and $\widetilde\ell_{\beta,\eta}[c]$ as the first derivative of $\widetilde\ell_\beta$ w.r.t. $\beta$ and $\eta$ (along the direction $c$), respectively. By setting $\epsilon=(\epsilon_\beta',0,\ldots,0)'$ and $\epsilon=(0,e)'=(0,\Delta H,b_1,\ldots,b_d)'$, respectively, some calculations reveal that
\begin{eqnarray}
E\left\{\widetilde\ell_{\beta,\beta}(X;\alpha_0^\ast)\epsilon_\beta\right\}+E\left\{\widetilde\ell_{\beta}(X;\alpha_0^\ast)\dot
\ell_{\beta}'(X;\alpha_0^\ast)\epsilon_\beta\right\}&=&0,\label{inter04}\\
E\left\{\widetilde\ell_{\beta,\eta}[e](X;\alpha_0^\ast)\right\}+E\left\{\widetilde\ell_{\beta}(X;\alpha_0^\ast)\dot\ell_{\eta}'[e](X;\alpha_0^\ast)
\right\}&=&0\label{inter05}
\end{eqnarray}
based on (\ref{inter03}). By considering the orthogonal property of $\widetilde\ell_{\beta_0}$ and the above reparametrization, we obtain the following two useful facts:
\begin{eqnarray}
&&\widetilde I_0=-E\left\{\widetilde\ell_{\beta,\beta}(X;\alpha_0^\ast)\right\},\label{effinf}\\
&&E\left\{\widetilde\ell_{\beta,\eta}[e](X;\alpha_0^\ast)\right\}=0\label{inter06}
\end{eqnarray}
based on (\ref{inter04}) and (\ref{inter05}).

Define $\widetilde\ell_{\beta,\alpha^\ast,\alpha^\ast}[h_1,h_2](X;\alpha^\ast)$ as the second order Fr\'{e}chet derivative of $\widetilde\ell_\beta$ w.r.t. $\alpha^\ast$ along the direction $[h_1,h_2]$ at the point $\alpha^\ast$. The same notation rule applies to $\dot\ell_{\beta,\alpha^\ast,\alpha^\ast}[h_1,h_2](X;\alpha^\ast)$ and $\dot\ell_{\eta,\alpha^\ast,\alpha^\ast}[h_1,h_2,h_3](X;\alpha^\ast)$. Now we are ready to express the Taylor expansion as follows.
\begin{eqnarray*}
&&E[\widetilde\ell_{\beta}(X;\alpha)-\widetilde\ell_{\beta}(X;\alpha_0)]\\
&=&E[\widetilde\ell_{\beta}(X;\alpha^\ast)-\widetilde\ell_{\beta}(X;\alpha_0^\ast)]\\
&=&E\left\{\widetilde\ell_{\beta,\beta}(X;\alpha_0^\ast)\right\}(\beta-\beta_0)+E\left\{\widetilde\ell_{\beta,\eta}[\eta-\eta_0](X;\alpha_0^\ast)\right\}\\
&&+\frac{1}{2}E\left\{\widetilde\ell_{\beta,\alpha^\ast,\alpha^\ast}[\Delta\alpha^\ast,\Delta\alpha^\ast](X;\widetilde\alpha^\ast)\right\}\\
&=&-\widetilde I_0(\beta-\beta_0)\\&&+\frac{1}{2}E\left\{\dot\ell_{\beta,\alpha^\ast,\alpha^\ast}[\Delta\alpha^\ast, \Delta\alpha^\ast](X;\widetilde\alpha^\ast)-\dot\ell_{\eta,\alpha^\ast,\alpha^\ast}[\bar c^\dag, \Delta\alpha^\ast,\Delta\alpha^\ast](X;\widetilde\alpha^\ast)\right\},
\end{eqnarray*}
where  $\Delta\alpha^\ast=\alpha^\ast-\alpha_0^\ast$ and $\widetilde\alpha^\ast$ lies between $\alpha^\ast$ and $\alpha_0^\ast$. The last equation in the above follows from (\ref{effinf}) \& (\ref{inter06}). Now we only need to show that the second term in the last equation is of the order
$$o(\|\beta-\beta_0\|)+o(n^{-1/2}).$$ Let $\Delta H=H-H_0$ and $\Delta h_j=h_j-h_{j0}$. After some algebra, we obtain
\begin{eqnarray*}
&&\dot\ell_{\beta,\alpha^\ast,\alpha^\ast}[\Delta\alpha^\ast,\Delta\alpha^\ast](X;\widetilde\alpha^\ast)\\
&=&Z\ddot Q_{\widetilde\theta}\left[Z'(\beta-\beta_0)+\Delta H(V)+\sum_{j=1}^d \Delta h_{j}(W_j)\right]^2,\\
&&\dot\ell_{\eta,\alpha^\ast,\alpha^\ast}[\bar c^\dag, \Delta\alpha^\ast,\Delta\alpha^\ast](X;\widetilde\alpha^\ast)\\
&=&\left[\int_{l_v}^V\bar a^\dag(s)dH(s)+\sum_{j=1}^d\bar b_j^\dag(W_j)\right]\ddot Q_{\widetilde\theta}\left[Z'(\beta-\beta_0)+\Delta H(V)+\sum_{j=1}^d \Delta h_j(W_j)\right]^2\\
&&+2\left[\int_{l_v}^V\bar a^\dag(s)d\Delta H(s)\right]\dot Q_{\widetilde\theta}\left[Z'(\beta-\beta_0)+ \Delta H(V)+\sum_{j=1}^d \Delta h_j(W_j)\right],
\end{eqnarray*}
where $\widetilde \theta$ lies between $\theta$ and $\theta_0$. Considering the assumption that $d(\alpha,\alpha_0)\leq C_1n^{-r/(2r+1)}$ and the previously shown result that $a_k^\dag$ and $b_{jk}^\dag$ are both uniformly bounded, we can verify Assumption B3 based on the above expressions. This completes the proof of Theorem~\ref{asymnor}. $\Box$

\subsection*{Proof of Theorem~\ref{infocon}}
For simplicity, we write $S_k(X;\alpha_0, w_k)$ and $S_k(X;\widehat\alpha,w_k)$ as $S_k^0[w_k]$ and
$\widehat S_k[w_k]$, respectively.
Based on the definitions of $\widetilde I_0$ and (\ref{ihat2}), we know their $(k,k')$-th entry can be written as
\begin{eqnarray}
\widetilde I_0(k,k')&=&ES_k^0[w_k^\dag]S_{k'}^0[w_{k'}^\dag],\label{iokk}\\
\widehat I(k,k')&=&\mathbb{P}_n\widehat S_k[\widehat w_k^\dag]\widehat S_{k'}[\widehat w_{k'}^\dag],\label{ihkk}
\end{eqnarray}
where $w_k^\dag=(a_k^\dag, b_{1k}^\dag,\ldots,b_{dk}^\dag)$ and $\widehat w_k^\dag=((\gamma_{0k}^\dag)'\mathbf{B}_0,(\gamma_{1k}^\dag)'\mathbf{B}_1,\ldots,(\gamma_{dk}^\dag)'\mathbf{B}_d)$. It is easy to show that
\begin{eqnarray}
E\left[\sup_{\alpha\in\mathcal{N}_0,w_k\in\mathcal{W}_n}|S_k(X;\alpha,w_k)|^2\right]\leq\mbox{const.}<\infty\label{bddcon}
\end{eqnarray}
since $\mathcal{A}$ and $\mathcal{W}_k$ are both assumed to be compact. Note that (\ref{bddcon}) implies that $\{S_k(x;\alpha,w_k):\alpha\in\mathcal{N}_0, w_k\in\mathcal{W}_{n}\}$ is P-Glivenko-Cantelli. Then, we know that, uniformly over $w_{k},w_{k'}\in\mathcal{W}_{n}$,
\begin{eqnarray}
&&\mathbb{P}_n\widehat S_k[w_{k}]\widehat S_{k'}[w_{k'}]\nonumber\\
&=&E\widehat S_k[w_{k}]\widehat S_{k'}[w_{k'}]+o_P(1)\label{setp1}
\end{eqnarray}
by considering Corollary 9.27 of \cite{k08}. Uniformly over $w_{k},w_{k'}\in\mathcal{W}_{n}$, we have
\begin{eqnarray}
&&\left|E\widehat S_k[w_{k}]\widehat S_{k'}[w_{k'}]-ES_k^0[w_{k}]S_{k'}^0[w_{k'}]\right|\nonumber\\
&\leq&E\left|\widehat S_k[w_{k}](\widehat S_{k'}[w_{k'}]-S_{k'}^0[w_{k'}])\right|+E\left|S_{k'}^0[w_{k'}](\widehat S_k[w_{k}]-S_k^0[w_{k}])\right|\nonumber\\
&\leq&\|\widehat S_k^2[w_{k}]\|_2\|\widehat S_{k'}[w_{k'}]-S_{k'}^0[w_{k'}]\|_2+\|S_{k'}^0[w_{k'}]\|_2\|\widehat S_k[w_{k}]-S_k^0[w_{k}]\|_2\nonumber\\
\;\;\;\;\;\;\;\;&\leq&o_P(1),\label{step2}
\end{eqnarray}
where the last inequality follows from (\ref{m41}) (together with the consistency of $\widehat \alpha$) \& (\ref{bddcon}).
Combining (\ref{setp1}) and (\ref{step2}), we have obtained that
\begin{eqnarray}
\sup_{w_{k},w_{k'}\in\mathcal{W}_{n}}\left|\mathbb{P}_n\widehat S_k[w_k]\widehat S_{k'}[w_{k'}]-ES_k^0[w_k]S_{k'}^0[w_{k'}]\right|=o_P(1),\label{step3}
\end{eqnarray}
which implies that
\begin{eqnarray}
\widehat I(k,k')=ES_k^0[\widehat w_k^\dag]S_{k'}^0[\widehat w_{k'}^\dag]+o_P(1).\label{step5}
\end{eqnarray}

To finish the proof, we need to introduce $\widetilde w_k^\dag\equiv\arg\min_{w_k\in\mathcal{W}_{n}}E\{S_k^0[w_k]\}^2$ as a bridge. Now, it remains to show that
\begin{eqnarray}
ES_k^0[\widehat w_k^\dag]S_{k'}^0[\widehat w_{k'}^\dag]-ES_k^0[\widetilde w_k^\dag]S_{k'}^0[\widetilde w_{k'}^\dag]&=&o_P(1),\label{step1}\\
ES_k^0[\widetilde w_k^\dag]S_{k'}^0[\widetilde w_{k'}^\dag]-\widetilde I_0(k,k')&=&o(1).\label{step4}
\end{eqnarray}
We first consider (\ref{step1}). By similar analysis applied to (\ref{step2}), we know that (\ref{step1}) holds if $\|S_k^0[\widetilde w_k^\dag]-S_k^0[\widehat w_k^\dag]\|_2=o_P(1)$. Denote $M_n(w)$ and $M(w)$ as $\mathbb{P}_n\widehat S_k^2[w]$ and $\|S_k^0[w]\|_2^2$, respectively.
The definition of $\widetilde w_k^\dag$ further implies that
\begin{eqnarray*}
\|S_k^0[\widetilde w_k^\dag]-S_k^0[\widehat w_k^\dag]\|_2^2&=&\|S_k^0[\widehat w_k^\dag]\|_2^2-\|S_k^0[\widetilde w_k^\dag]\|_2^2,\\
&=&\mathbb{P}_n\widehat S_k^2[\widehat w_k^\dag]-\|S_k^0[\widetilde w_k^\dag]\|_2^2+o_p(1),\\
&=&M_n(\widehat w_k^\dag)-M(\widetilde w_k^\dag)+o_P(1),
\end{eqnarray*}
where the second equality follows from (\ref{step3}). By the definitions of $\widehat w_k^\dag$ and $\widetilde w_k^\dag$, we have
$$M_n(\widehat w_k^\dag)-M(\widehat w_k^\dag)\leq M_n(\widehat w_k^\dag)-M(\widetilde w_k^\dag)\leq M_n(\widetilde w_k^\dag)-M(\widetilde w_k^\dag).$$ Therefore, we conclude the proof of (\ref{step1}) by applying (\ref{step3}) to the above inequality. We next consider (\ref{step4}). Again, by the form of $\widetilde I_0(k,k')$ given in (\ref{iokk}) and similar analysis in (\ref{step1}), we only need to show $\|S_k^0[\widetilde w_k^\dag]-S_k^0[w_k^\dag]\|_2=o(1)$. By the definitions of $\widetilde w_{k}^\dag$ and $w_{k}^\dag$, we have
\begin{eqnarray*}
\|S_k^0[\widetilde w_k^\dag]-S_k^0[w_k^\dag]\|_2^2&=&\inf_{w_k\in \mathcal{W}_{n}}E\left[\dot\ell_{g_0}[a_k^\dag]-\dot\ell_{g_0}[a_k]+\sum_{j=1}^d(\dot\ell_{h_{j0}}[b_{jk}^\dag]-\dot\ell_{h_{j0}}[b_{jk}])\right]^2\\
&\aplt&\inf_{w_k\in \mathcal{W}_{n}}\left\{\|\dot\ell_{g_0}[a_k^\dag]-\dot\ell_{g_0}[a_k]\|_2^2+\sum_{j=1}^{d}\|\dot\ell_{h_{j0}}[b_{jk}^\dag]-
\dot\ell_{h_{j0}}[b_{jk}]\|_2^2\right\}\\
&\aplt&\inf_{a_k\in\mathcal{G}_n}\|\dot\ell_{g_0}[a_k^\dag]-\dot\ell_{g_0}[a_k]\|_2^2+\sum_{j=1}^d\inf_{b_{jk}\in \mathcal{H}_{jn}}\|\dot\ell_{h_{j0}}[b_{jk}^\dag]-\dot\ell_{h_{j0}}[b_{jk}]\|_2^2\\
&\aplt&\inf_{a_k\in\mathcal{G}_n}\|a_k^\dag-a_k\|_\infty^2+\sum_{j=1}^d\left\{\inf_{b_{jk}\in \mathcal{H}_{jn}}\|b_{jk}^\dag-b_{jk}\|_\infty^2\right\},
\end{eqnarray*}
where the last inequality trivially follows from the form of $\dot\ell_g[a]$ and $\dot\ell_{h_j}[b_j]$. According to the analysis in the proof of Theorem~\ref{asymnor}, we know that $a_k^\dag\in H_{\widetilde c_0}^{r_0}[l_v,u_v]$ and $b_{jk}^\dag\in H_{\widetilde c_j}^{r_j}[0,1]$. Thus, we have $\|S_k^0[\widetilde w_k^\dag]-S_k^0[w_k^\dag]\|_2\rightarrow 0$ based on the last inequality in the above. This completes the whole proof. $\Box$

\bibliographystyle{biometrika}

\end{document}